\documentclass[reqno]{amsart}

\usepackage{amssymb}
\usepackage{amscd}
\usepackage{color}
\usepackage{amsfonts}
\usepackage{eucal}
\usepackage{latexsym}
  \RequirePackage{amsbsy}

\newcommand{\GCD}{\text{\rm GCD}}

\newcommand{\Prin}{\text{\rm Prin}}

\newcommand{\Inv}{\text{\rm Inv}}

\newcommand{\f}{\boldsymbol{f}}
\newcommand{\F}{\boldsymbol{F}}
\newcommand{\astf}{\ast_{\!{_f}}}
  
\newtheorem{theorem}{Theorem}[section]

\newtheorem{corollary}[theorem]{Corollary}

\newtheorem{example}[theorem]{Example}
\newtheorem{remark}[theorem]{Remark}
\newtheorem{proposition}[theorem]{Proposition}

\title[On $v$--domains and star operations]
{On $v$--domains and star operations}

\author[D.D. Anderson, D.F. Anderson, M. Fontana, and M. Zafrullah]{ D.D. Anderson, David F. Anderson, \\ Marco Fontana and  Muhammad Zafrullah}
\address{D.D.A.: Department of Mathematics, University of
Iowa, Iowa City, IA 52242, USA.}

\thanks{\it Acknowledgments. \rm  We would like to thank Franz Halter-Koch for several helpful comments on a previous version of the present paper and for suggesting   lines for further investigation and extension to the general setting of ideal  systems on monoids.  \\
$\mbox{ \quad}$ During the preparation of this paper, the third named author
was partially supported by  a research grant PRIN-MiUR}
\email{dan-anderson@uiowa.edu}
\address{ D.F.A.: Department of Mathematics, University of
Tennessee at Knoxville, Knoxville, TN 37996-1300, USA.}
\email{anderson@math.utk.edu}
\address{M.F.: \ Dipartimento di Matematica, Universit\`a degli Studi
``Roma Tre'', 00146 Rome, Italy.}
\email{fontana@mat.uniroma3.it }
\address{M.Z.: Department of Mathematics, Idaho State University, Pocatello, ID 83209, USA.}
\email{mzafrullah@usa.net}

\date{\today }

\subjclass[2000]{ 13A15, 13F05, 13G05.  }

\keywords{ Pr\"ufer domain, GCD domain, star operation, $v$--domain,  $\ast$--Pr\"ufer domain, completely
integrally closed domain,  $\ast$--invertible ideal.  }

\newtheorem{problem}[theorem]{Problem}

\begin{document}

\maketitle

  \begin{abstract} Let $\ast $ be a star operation on an integral domain $D$.
Let $\f(D)$ be the set of all nonzero finitely generated fractional ideals of $D$. Call $D$ a    $\ast$--Pr\"ufer (respectively, $(\ast, v)$--Pr\"ufer) domain         if $(FF^{-1})^{\ast}=D$ (respectively, $(F^vF^{-1})^{\ast}=D$)  for  all $F\in \f(D)$.  We establish  that   $\ast$--Pr\"ufer domains   (and  $(\ast, v)$--Pr\"ufer domains)   for various  star operations   $\ast $ span a major
portion of the known generalizations of Pr\"{u}fer domains  inside the class of $v$--domains.   We also  use Theorem 6.6 of  the  Larsen and   McCarthy book  [Multiplicative Theory of Ideals, Academic Press, New York--London, 1971],  which gives several equivalent conditions for an integral domain to be a Pr\"ufer domain,  as a  model,  and we show
which statements of that theorem on Pr\"ufer domains can be generalized in a natural way and  proved for   $\ast$--Pr\"ufer domains,   and
which cannot be.  We  also show that  in a   $\ast $--Pr\"ufer domain,  each pair of     $\ast $-invertible $\ast $-ideals   admits  a GCD in the set of $\ast $-invertible $\ast $-ideals, obtaining a remarkable generalization of a property holding for the ``classical'' class of Pr\"ufer $v$--multiplication domains.   We also
 link $D$ being   $\ast $--Pr\"ufer (or  $(\ast, v)$--Pr\"ufer)  with the  group Inv$^{\ast }(D)$ of $\ast $-invertible $\ast $-ideals (under $\ast$-multiplication)  being  lattice-ordered.   
  \rm

 \end{abstract}

 \bigskip

 The so called $v$--domains (i.e., the integral domains  such that every nonzero finitely generated fractional ideal is $v$--invertible) include several distinguished classes of Pr\"ufer-like domains, but not
much seems to be known about them. (For a brief history of $v$-domains and   an annotated   list of references on this important class of domains,  see \cite{za-v-domains}).  The aim of this article is to prove
new properties of $v$--domains in their most general form, using star operations, and to   give  a  unifying pattern in  this body of results.  As a   consequence,   after specializing the star operation to some relevant cases, we also   obtain several  properties already known for various classes of Pr\"ufer-like domains, providing a clear indication how these properties and classes of domains are related to one another.

 Let $D$ be an integral domain with quotient field  $K$,   and let $\F(D)$ denote
the set of nonzero fractional ideals of $D$.   Also,   let $\f(D):=\{A\in \F(D)\mid A$
is finitely generated$\}$ and $\F^v(D) :=\{ A \in \F(D) \mid A = A^v\}$  the set of fractional divisorial ideals,    where   $A^v := \left(A^{-1}\right)^{-1}$.   Let $\ast $ be a star operation on $D$.    (For a   review  of star  operations,   the reader may consult  Gilmer \cite[Sections 32 and 34]{G}  or Halter-Koch \cite{hkbook} for a  general approach in the language of   ideal systems on monoids.)  Call $
A\in \F(D)$ \it $\ast$-invertible \rm  if $(AA^{-1})^{\ast }=D$.  We say that $D$ is a \it 
completely integrally closed  domain \rm  (for short, CICD) if  $ D = \widetilde{D} :=
 \{ x\in K \mid $  there exists $ 0 \neq d\in D $  such that $dx^{n}\in D$ 
for all  integers $n\geq 1 \}$.  
It is well known that $D$ is a 
completely integrally closed domain    if and only if   $
D = (AA^{-1})^{v} \  (=(A^{v}A^{-1})^{v})$   for all $A\in \F(D)$ \cite[Theorem 34.3]{G}; in particular, a CICD is a $v$--domain. 

These days it is customary to take a concept defined or
characterized using the standard $v$-operation and to ask for domains that
are characterized or defined by replacing the $v$-operation by a general
star operation.  It appears that  in the case of CICD's there are at least two
 star operation analogues.  
Let $\ast $ be a general star operation on $D$.
Call $D$ a \it $\ast $--completely integrally closed domain \rm (for short,  $\ast $--CICD) if $(AA^{-1})^{\ast }=D$ for all   $A\in \F(D)$,   and
call $D$ a  \it $(\ast, v) $--completely integrally closed domain \rm (for short,  $(\ast, v) $--CICD)  if $(A^{v}A^{-1})^{\ast }=D$ for all $A\in \F(D)$.    Clearly  a $\ast$--CICD,  or a $(\ast, v)$--CICD, is a $v$--domain, since $(A^\ast)^v = A^v = (A^v)^\ast$ for all $A\in \F(D)$. 
 Moreover,  
a  $\ast$--CICD\  is a  $(\ast, v)$--CICD,   but not conversely. 

 Let  $\ast_1, \ast_2$ be two star operations   on $D$. Recall that  $\ast_1 \leq  \ast_2$ if $A^{\ast_1} \subseteq A^{\ast_2}$  for all $A\in \F(D)$. For instance, we have $w \leq t \leq v$, where the \it $t$--operation \rm (respectively, {\it{$w$--operation}}) is defined by setting $A^t := \bigcup \{F^v \mid F \subseteq A,\  F \in \f(D) \}$ (respectively,   $A^w := \bigcap \{ AD_P \mid P \mbox{ is a maximal $t$--ideal of $D$} \}$) for all $A \in \F(D)$; cf. for instance \cite[Theorem 34.1]{G}, \cite{WM}, \cite[Section 2]{AC}, and \cite[Sections 3 and 4]{FH}.  

Clearly, if $\ast_1, \ast_2$ are two star operations   on $D$  with   $\ast_1 \leq  \ast_2$, then  a $\ast_1 $--CICD (respectively, $(\ast_1, v) $--CICD)  $D$  is a  $\ast_2 $--CICD (respectively, $(\ast_2, v) $--CICD). 
Recall  that   $\ast \leq v$  for  every  star operation $\ast$  on $D$  \cite[Theorem 34.1]{G},  and  hence $(AA^{-1})^{\ast }=D$ implies $(AA^{-1})^{v}=D$ for $A \in \F(D)$, i.e.,  a $\ast $-invertible
ideal is $v$-invertible.    Therefore a $(\ast, v) $--CICD (and, in particular, a $\ast $--CICD) is a $v$--CICD (= CICD).   However,  since $A\in \F(D)$ being $\ast $-invertible implies $A^{\ast
}=A^{v}$ \cite[page 433]{Z-inv},    a  distinction between $\ast $--CICD's and $(\ast, v)$--CICD's appears highly unlikely. But, as we shall see, there is a marked
distinction   between them  in several cases.

In a preliminary part  of this   paper,  we discuss the motivations and the advantages for studying  star operation analogues of CICD's,    we  give some characterizations of $(\ast, v)$--CICD's,   we  
 give interpretations of $(\ast , v)$--CICD's for different star operations $\ast $,
  we    compare them with $\ast$--CICD's, and
  we  review results known for both.

 Having dealt with this topic of immediate
interest in Section 1,  
  we investigate   in Section 2   the main theme of this paper studying   a  ``star operation version''   of $v$--domains.
 We call $D$ a  \it   $\ast$--Pr\"ufer domain \rm       if every nonzero finitely generated  ideal of $D$ is 
$\ast $-invertible  (i.e.,  $(FF^{-1})^{\ast}=D$ for  all $F\in \f(D)$),   and we call $D$ a  \it  $(\ast, v)$--Pr\"ufer domain \rm     if $F^{v}$ is $\ast $-invertible   (i.e., $(F^vF^{-1})^{\ast}=D$)  for  all  $F \in \f(D)$.   Clearly, if $\ast_1, \ast_2$ are two star operations  on $D$  with   $\ast_1 \leq  \ast_2$, then a  $\ast_1 $--Pr\"ufer domain (respectively, $(\ast_1, v) $--Pr\"ufer  domain)  $D$  is a $\ast_2 $--Pr\"ufer  domain (respectively, $(\ast_2, v) $--Pr\"ufer  domain).   Clearly  a $\ast$--Pr\"ufer  domain,  or a $(\ast, v)$--Pr\"ufer  domain, is a $v$--domain. Moreover,  
  a   $\ast$--Pr\"ufer domain  is a $(\ast, v)$--Pr\"ufer domain,     but not conversely. 
 
 These domains have been  partially    studied in \cite{AMZ} as
special cases    of rather general results \cite[Theorem 4.1 and Corollary 4.3]{AMZ}.  Since  the proofs provided in \cite{AMZ}  
were sort of dismissive, we provide here direct proofs of some  of   the relevant  results stated in
\cite{AMZ} and we prove some more. We establish in this section that   $\ast$--Pr\"ufer domains   (and  $(\ast, v)$--Pr\"ufer domains)   for various  star operations   $\ast $ span a major
portion of the known generalizations of Pr\"{u}fer domains  inside the class of $v$--domains.   For example, for   $\ast = d$   (i.e., 
the identity star operation), we get a $d$--Pr\"ufer
domain which is precisely a Pr\"ufer   domain; for $\ast =t$,    we get a $t$--Pr\"ufer
domain which is precisely a Pr\"ufer $v$--multiplication domain (or a  P$v$MD);   and of
course for   $\ast =v$,   we get the usual $v$--domain.     In this  section,  
we also  use Theorem 6.6 of Larsen and   McCarthy  \cite{LM}, which gives several equivalent conditions for an integral domain to be a Pr\"ufer domain,  as a  model,  and we show
which statements of that theorem on Pr\"ufer domains can be generalized in a natural way and  proved for   $\ast$--Pr\"ufer domains   and
which cannot be.   In particular, we show that $D$ is  a $\ast $--Pr\"ufer domain  if and only if $((A\cap
B)(A+B))^{\ast }=(AB)^{\ast }$ for all $A,B\in \F(D)$. This type of  result is  known
for  Pr\"ufer domains \cite[Theorem 25.2]{G} and for P$v$MD's  \cite[Theorem 5]{Gr},  but is definitely not
known for $v$-domains.

   The last part of the paper deals with a general form of GCD for $\ast$--Pr\"ufer domains  (in particular, for $v$--domains) and connections with   lattice-ordered   abelian groups. The key  fact is that   an integral domain $D$ is a 
$\ast$--Pr\"ufer domain if and only if $A+B$ is $\ast$-invertible for all $
\ast $-invertible $A, B\in \F(D)$. \
Recall that  $D$ is  a \it GCD domain \rm  if for  all     $
x,y\in D^{\times }:=D\backslash \{0\}$,   we have $\GCD(x,y)\in D$.  \
Now a B\'ezout domain $D$ (e.g., a PID) is slightly
more than a GCD domain in that for  all   nonzero ideals $aD$  and  $bD$,   we
have a unique ideal $dD$ with    $aD+bD=dD$,    where   $d$ is a GCD of  $a$ and $b$. 
 Moreover,  note
 that  $aD$ and $ bD$   are invertible ideals and that in a Pr\"{u}fer domain nonzero finitely generated ideals are  invertible.   If we
regard, for every pair of   integral    invertible ideals $A$ and $B$ of a Pr\"ufer domain, the invertible ideal $
C:=A+B$ as the GCD of $A$ and  $B$,  then   we find that   $A,B \subseteq C$.   Hence, 
 $A_{1}:= AC^{-1} \subseteq D$ and   $B_{1}:= BC^{-1} \subseteq D$,    and so $A=A_{1}C$ and  $B=B_{1}C$,   where $A_{1}+B_{1}=D.$ Thus, in a Pr\"{u}fer  domain,   each pair of  integral  invertible ideals
has  a   GCD of sorts. In  \cite[Section 1]{BGZ},     the above observations were used to show that
in    $t$--Pr\"ufer domains  (= Pr\"ufer $v$--multiplication domains),    each  pair of   integral  $t$-invertible $t$-ideals has a
GCD of sorts, generalizing to this setting some aspects of the GCD theory of B\'ezout domains.  
 In the general context of  $\ast $--Pr\"ufer domains, we  show   that each pair of  integral  $\ast $-invertible $\ast $-ideals  has   a GCD of sorts.   This result is a slightly bigger jump than the $t$--Pr\"ufer  domain case in that, in a  $\ast $--Pr\"ufer domain,  a $\ast $-invertible $\ast $-ideal may not be of finite type.

\section{Star Completely Integrally Closed Domains}

\bigskip 
  Recall that a ring $R$ is a \it multiplication ring \rm   if   for  all  ideals $A$ and $ B$   of $R$ with $A \subseteq B$, there exists an ideal $C$ of $R$ such that $A = BC$ (cf. for  instance,   \cite[Definition 9.12, page 209]{LM}). Clearly, a Dedekind domain is a  multiplication    ring,   and more precisely, for an integral domain the notions of Dedekind domain and  multiplication  ring coincide \cite[Theorem 9.13]{LM}.

Given a star operation $\ast$ on an integral domain  $D$,   it is natural to call $D$ a \it  $\ast$--multiplication domain \rm   (respectively, \it $(\ast, v)$--multiplication domain\rm ) if for
 all  $A,B\in $ $\F(D)$ with $A^\ast \subseteq B^\ast$  (respectively,  $A^\ast \subseteq B^v$),   there exists $C \in \F(D)$ such that $A^\ast = (BC)^\ast$  (respectively,  $A^\ast = (B^vC)^\ast$).  Note that star multiplication  domains,  and in particular,  divisorial multiplication domains were recently investigated in relation  to   Gabriel topologies by J. Escoriza and B. Torrecillas \cite{ET}.
 
    As usual, for all $A,B\in \F(D)$, we denote  the fractional ideal  $(A:_K B) := \{x \in K \mid xB \subseteq A\}$  by 
$(A:B)$  and  the ideal $(A:B) \cap D$ by $(A:_D B)$.

\begin{proposition} \label{prop:1.0}
Let $\ast $ be a star operation on  an integral domain  $D.$ Then
\begin{enumerate}
\item[(a)]  $D$
is a $\ast$--CICD if and only if $D$ is a $\ast$--multiplication  domain, \, and  
\item[(b)]  $D$
is a $(\ast, v)$--CICD if and only if $D$ is a $(\ast, v)$--multiplication  domain.  
\end{enumerate}
\end{proposition}
\begin{proof} (a)  If $D$ is   a  $\ast$--CICD and  $A^\ast \subseteq B^\ast$,   then $C:= B^{-1}A \in \F(D)$  satisfies   $(BC)^\ast =(BB^{-1}A)^\ast =  ((BB^{-1})^\ast A)^\ast =(DA)^\ast = A^\ast$.

Conversely, for each  $A \in \F(D)$,  let   $0\neq a \in A$,  and so $aD \subseteq A^\ast$. By  assumption,   there exists $C \in \F(D)$ such that $(AC)^\ast = aD$,   i.e.,  $(Aa^{-1}C)^\ast = D$. Note that ~$B:= a^{-1}C \in \F(D)$ and $a^{-1}C  \subseteq (D:A)$. Therefore we conclude that $D = (Aa^{-1}C)^\ast \subseteq (AA^{-1})^\ast \subseteq D$,  i.e.,   $(AA^{-1})^\ast  = D$.

(b)   If $D$ is  a  $(\ast, v)$--CICD and  $A^\ast \subseteq B^v$,   then $C:= B^{-1}A \in \F(D)$  satisfies   $(B^vC)^\ast =(B^vB^{-1}A)^\ast =  ((B^vB^{-1})^\ast A)^\ast =(DA)^\ast = A^\ast$.

Conversely, for each  $A \in \F(D)$,  let $0\neq a \in A$,   and so $aD \subseteq A^v$.   By assumption,   there exists $C \in \F(D)$ such that $(A^vC)^\ast = aD$,   i.e.,  $(A^va^{-1}C)^\ast = D$. Note that $B:= a^{-1}C \in \F(D)$ and $a^{-1}C  \subseteq (D:A^v) =(D:A)$. Therefore we conclude that $D = (A^va^{-1}C)^\ast \subseteq (A^vA^{-1})^\ast \subseteq D$,   i.e.,   $(A^vA^{-1})^\ast  = D$.
\end{proof}

\begin{proposition} \label{prop:1.1}
Let $\ast $ be a star operation on  an integral domain  $D$. Then $D$
is a $(\ast, v)$--CICD if and only if $(AB)^{-1}=(A^{-1}B^{-1})^{\ast }$ for
  all   $A,B\in \F(D).$
\end{proposition}

\begin{proof} Suppose that $D$ is a $(\ast, v)$--CICD and consider $A,B\in \F(D)$.
Then $D=((AB)^{v}(AB)^{-1})^{\ast }$.  Multiplying both sides of the above
equation by $A^{-1}B^{-1}$ and applying $\ast $,  we get:
$$
\begin{array}{rl}
(A^{-1}B^{-1})^{\ast
}= & \hskip -5pt (A^{-1}B^{-1} (AB)^{v}(AB)^{-1})^{\ast }   \supseteq 
(A^{-1}B^{-1}A^{v}B^{v}(AB)^{-1})^{\ast }  \\
= &\hskip -5pt  ((A^{-1}A^{v})(B^{-1}B^{v})(AB)^{-1})^\ast \supseteq  (A^{-1}A^{v})^{\ast
}(B^{-1}B^{v})^{\ast }((AB)^{-1})^{\ast } \\ =& \hskip -5pt   ((AB)^{-1})^{\ast }=(AB)^{-1}\,.
\end{array}
$$

For the reverse   inclusion,  note that $A^{-1}B^{-1}\subseteq (AB)^{-1}$,   and so 
$(A^{-1}B^{-1})^{\ast }\subseteq  ((AB)^{-1})^{\ast }= (AB)^{-1}$. 

Conversely, if $(AB)^{-1}=(A^{-1}B^{-1})^{\ast }$ for all $A,B  \in \F(D)$, 
then in particular $(A^{-1}A^{v})^{-1} = (A^{v}A^{-1})^{\ast }$ for all $A\in
\F(D)$.  Now, as $A^{v}A^{-1}=A^{-1}A^{v}\subseteq D$,  we have $D\subseteq
(A^{-1}A^{v})^{-1}\ = (A^{v}A^{-1})^{\ast } \subseteq D$.   Thus $(A^{v}A^{-1})^{\ast } = D$ for all  $A\in
\F(D)$;  so $D$ is a $(\ast, v)$--CICD. 
\end{proof}


\begin{proposition} \label{cor:1.2}  Let $\ast $ be a star operation on  an integral domain  $D$. 
Then
 \begin{enumerate}
\item[(a)] If $D$ is a
$\ast $--CICD, then $D$  is a $(\ast, v)$--CICD.

\item[(b)]  If $D$ is a $(\ast, v)$--CICD, then $D$ is a completely integrally  closed domain.
\end{enumerate}
 \end{proposition}

 \begin{proof}  From the definition and from the fact that $(A^v)^{-1} =A^{-1}$ for all $A \in \F(D)$, it follows immediately that  a $\ast $--CICD is a $(\ast, v)$--CICD.  Furthermore, if  we have  $(A^{v}A^{-1})^{\ast  } = D$ for all $A\in \F(D)$,
  then 
 $(A^{v}A^{-1})^{v}=D$,  and so  $(AA^{-1})^{v}=(A^{v}A^{-1})^{v}=D$, i.e., $D$ is a completely  integrally  closed domain.
\end{proof}

These results are simple and straightforward, but their value is in the
interpretation of the $(\ast, v)$--CICD for different star operations $\ast$.
We shall give examples of $(\ast, v)$--CICD's that are not 
$\ast $--CICD's for
the same $\ast $. 
Most of our examples come from  
\cite{AMZ},  which provides a lot of  quotient-based  characterizations of $\ast $--CICD's and of $(\ast, v)$--CICD's. Since the method of proof in \cite{AMZ}  was
somewhat  involved,   we include direct proofs of  these characterizations here.

\begin{proposition} \label{*CICD} \rm \cite[Corollary 3.4]{AMZ}  \it Let $\ast $ be a star operation on an integral domain  $
D.$ Then the following conditions are   equivalent. 
\begin{enumerate}

\item[(i)] $(A:B)^{\ast }=(AB^{-1})^{\ast }$ for all $A,B\in \F(D)$.

\item[(ii)]  $(A:B^{-1})^{\ast }=(AB)^{\ast }$ for all $A,B\in \F(D)$.

\item[(iii)]  $(A^{\ast }:B) =(AB^{-1})^{\ast }$ for all $A,B\in \F(D)$.

\item[(iv)]  $(A^{\ast }:B^{-1})=(AB)^{\ast }$ for all $A,B\in \F(D)$.

\item[(v)]  $D$ is a $\ast $--CICD.

\item[(vi)]  $D$ is a CICD and $A^{\ast }=A^{v}$ for  all   $A\in \F(D)$.

\item[(vii)]  $(A^{v}:B^{-1})=(A^{v}B)^{\ast }$\ for all $A,B\in \F(D)$.
\end{enumerate}
\end{proposition}

\begin{proof}  Let us note that (v)$\Leftrightarrow $(vi) is well known and it is
the only part of the proof directly given in \cite[Corollary 3.4 and Proposition 3.2]{AMZ}. For the
 rest,   we use the following plan: (i)$\Rightarrow $(iii)$\Rightarrow $(v)$
\Rightarrow $(i) and (ii)$\Rightarrow $(iv)$\Rightarrow $(vii)$\Rightarrow 
$(v)$\Rightarrow $(ii).

(i)$\Rightarrow$(iii). Replace  $A$ by $A^{\ast }$  in (i)   to get   $(A^{\ast
}:B)^{\ast }=(A^{\ast }B^{-1})^{\ast }$,   and note that $(A^{\ast }:B)$ is a $
\ast $-ideal and that $(A^{\ast }B^{-1})^{\ast }=(AB^{-1})^{\ast }$.   So   $
(A^{\ast }:B)=(AB^{-1})^{\ast }$ for all $A,B\in \F(D).$

(iii)$\Rightarrow $(v). Set $B=A^{\ast }$ in (iii) to get $(A^{\ast }:A^{\ast
})=(A(A^\ast)^{-1})^{\ast }= (AA^{-1})^{\ast }$.   Noting that  $D\subseteq (A^{\ast }:A^{\ast
})=(AA^{-1})^{\ast }\subseteq D$,   we have $(AA^{-1})^{\ast }=D$ for all $A\in
 \F(D)$. 

(v)$\Rightarrow $(i). Note that  $AB^{-1}\subseteq (A:B)$;   so $
(AB^{-1})^{\ast }\subseteq (A:B)^{\ast }$.   For the reverse  inclusion,   let $
x\in (A:B)^{\ast }$.    Then   $xB\subseteq (A:B)^{\ast }B$,   and so   $xB^{\ast
}\subseteq ((A:B)^{\ast }B)^{\ast }=((A:B)B)^{\ast }\subseteq A^{\ast }.$
This gives $xB^{\ast }\subseteq A^{\ast }$. Multiplying by $B^{-1}$ on both
 sides   and applying $\ast $,   we have  $(xB^{\ast }B^{-1})^\ast \subseteq (A^{\ast }B^{-1})^\ast = (AB^{-1})^\ast  $.  Invoking (v),  we get $x\in (AB^{-1})^{\ast }.$

(ii)$\Rightarrow $(iv). Same as (i)$\Rightarrow $(iii).

(iv)$\Rightarrow $(vii). Replace $A$ by $A^{v}$ in (iv) to get  $((A^{v})^{\ast
}:B^{-1})=(A^{v}B)^{\ast }$,   and note that $(A^{v})^{\ast }=A^{v}$.

(vii)$\Rightarrow $(v). Set $A=B^{-1}$ in (vii) to get  $
(B^{-1}:B^{-1}) =(B^{-1}B)^{\ast }$,  and proceed as in the proof of (iii)$
\Rightarrow$(v) in order to get $(B^{-1}B)^{\ast }= D$.

(v)$\Rightarrow $(ii). The proof is more or less similar to the proof of
(v)$\Rightarrow $(i). More precisely,   $ AB \subseteq (A:B^{-1})$,   and so 
$ (AB)^\ast  \subseteq (A:B^{-1})^\ast$.
 Conversely,   let $
x\in (A:B^{-1})^\ast$.    Then   $xB^{-1}\subseteq (A:B^{-1})^{\ast }B^{-1}$,   and so $x(B^{-1})^{\ast
}\subseteq ((A:B^{-1})^{\ast }B^{-1})^{\ast }=((A:B^{-1})B^{-1})^{\ast }\subseteq A^{\ast }$.
This gives $x(B^{-1})^{\ast }\subseteq A^{\ast }$. Multiplying by $B$ on both  
 sides     and applying $\ast $,   we have  $(x(B^{-1})^{\ast }B)^\ast \subseteq (A^{\ast }B)^\ast = (AB)^\ast  $.  Invoking (v),  we get $x\in (AB)^{\ast }$.  \end{proof}

For the $(\ast, v)$--CICD  case,   we have the following set of  quotient-based  characte\-ri\-zations.

\begin{proposition}  \label{*v-cicd} \rm \cite[Corollary 3.5]{AMZ} \it   Let $\ast $ be a star operation
 on  an integral domain $D$.  Then the following  conditions are  equivalent. 
\begin{enumerate}

\item[(i)]  $D$ \ is a $(\ast, v)$--CICD.

\item[(ii)]  $(A^{v}:B)=(A^{v}B^{-1})^{\ast }$ for all $A,B\in \F(D)$. 

\item[(iii)]   $(A^{v}:B^{-1})=(A^{v}B^{v})^{\ast }$ for all $A,B\in \F(D)$.

\item[(iv)]   $(A:B)^{v}=(A^{v}B^{-1})^{\ast }$ for all $A,B\in \F(D)$.

\item[(v)]  $(A:B^{-1})^{\ast }=(AB^{v})^{\ast }$ for all $A,B\in \F(D)$.

\item[(vi)]  $(A^{\ast }:B^{-1})=(AB^{v})^{\ast }$ for all $A,B\in \F(D)$.
\end{enumerate}
\end{proposition}

\begin{proof}  (i)$\Rightarrow $(ii). Obviously  $A^{v}B^{-1}\subseteq (A^{v}:B)$,   and
so $(A^{v}B^{-1})^{\ast }\subseteq (A^{v}:B)^\ast = (A^{v}:B)$.  For the reverse  inclusion,   let $
x \in (A^{v}:B)$.  Then  $xB\subseteq A^{v}$,   and so $xB^{v}\subseteq A^{v}$.
Multiplying the last equation by $B^{-1}$ and applying $\ast $,  we get $
(xB^{v}B^{-1})^{\ast }\subseteq (A^{v}B^{-1})^{\ast }$.  Invoking  (i),   we have 
$x\in (A^{v}B^{-1})^{\ast }$,   and from this follows $(A^{v}:B)\subseteq
(A^{v}B^{-1})^{\ast }.$

(ii)$\Rightarrow $(iii). Replace $B$ by $B^{-1}$ in $(A^{v}:B)=(A^{v}B^{-1})^{
\ast}$ for all $A,B\in \F(D).$

(iii)$\Rightarrow $(i). Set $B=A^{-1}$ in the equality  $(A^{v}:B^{-1})=(A^{v}B^{v})^{\ast } 
$ to get $(A^{v}:A^{v})=(A^{v}A^{-1})^{\ast }$.   But since $D\subseteq
(A^{v}:A^{v})=(A^{v}A^{-1})^{\ast }\subseteq D$,  we conclude that $(A^{v}A^{-1})^{\ast }=D$  for all $
A\in \F(D)$,   and so $D$ is a   $(\ast , v)$--CICD.  

(i)$\Rightarrow $(iv). Note that $AB^{-1}\subseteq (A:B)$. So $
(AB^{-1})^{v}\subseteq (A:B)^{v}$.   This gives  $A^{v}B^{-1}\subseteq
(AB^{-1})^{v}\subseteq (A:B)^{v}$,   and from this we conclude that $
(A^{v}B^{-1})^{\ast }\subseteq (A:B)^{v}$.    
For the reverse inclusion,  let $
x\in (A:B)^{v}$. Then  $xB\subseteq (A:B)^{v}B$,   and so $xB^{v}\subseteq
((A:B)^{v}B)^{v}=((A:B)B)^{v}\subseteq A^{v}$.
 This gives  $xB^{v}\subseteq
A^{v}$, and   on multiplying by $B^{-1}$ on both sides,  we get $
xB^{v}B^{-1}\subseteq A^{v}B^{-1}$.   Applying $\ast $ on both sides and
invoking  (i),   we conclude that $x\in (A^{v}B^{-1})^{\ast }$.  This  
establishes the reverse inclusion.

(iv)$\Rightarrow $(i). Set $B=A$ in $(A:B)^{v}=(A^{v}B^{-1})^{\ast }$ to
get $(A:A)^{v}=(A^{v}A^{-1})^{\ast }$.     But then $D\subseteq 
(A:A)^{v}=(A^{v}A^{-1})^{\ast }\subseteq D$ for  all   $A\in \F(D)$.  That is, $
(A^{v}A^{-1})^{\ast }=D$ for  all  $A\in \F(D)$,   and this is (i).

(v)$\Rightarrow $(vi).  This is obvious  once we replace $A$ by $A^{\ast }$ and note
that $(A^{\ast }:B^{-1})^{\ast }=(A^{\ast }:B^{-1}).$

(vi)$\Rightarrow $(i). Set $A=B^{-1}$ in (vi) to get $
(B^{-1}:B^{-1})=(B^{-1}B^{v})^{\ast }$,  which can be used to conclude that $
(B^{-1}B^{v})^{\ast }=D$ for all $B\in \F(D)$.

(i)$\Rightarrow $(v). Clearly   $AB^{v}\subseteq( A:B^{-1})$,   and so $
(AB^{v})^{\ast }\subseteq (A:B^{-1})^{\ast }$.    For the reverse  inclusion,   let 
$x\in (A:B^{-1})^{\ast }$. Then $xB^{-1}\subseteq (A:B^{-1})^{\ast }B^{-1}$.
So $xB^{-1}\subseteq ((A:B^{-1})^{\ast }B^{-1})^{\ast
}=((A:B^{-1})B^{-1})^{\ast }\subseteq A^{\ast }$. Multiplying both sides of $
xB^{-1}\subseteq A^{\ast }$ by $B^{v}$,  we have $xB^{-1}B^{v}\subseteq
A^{\ast }B^{v}$.  Applying $\ast $ and invoking (i),  we conclude that $x\in
(A^{\ast }B^{v})^{\ast }=(AB^{v})^{\ast }$.
\end{proof}

\begin{remark} \rm
It is easy to verify that   statement   (ii) of Proposition \ref {*v-cicd} can be equivalently stated  as in \cite[Corollary 3.5 (2)]{AMZ}:

\begin{enumerate} \it
\item[(ii$^\prime$)]  $(A^{v}:B^{v})=(A^{v}B^{-1})^{\ast }$ for all $A,B\in \F(D)$. 
\end{enumerate}
\rm
\end{remark}

\medskip

The next result provides  a    useful characterization of $\ast$-invertible fractional ideals and sheds new light on Proposition \ref{*CICD}.

\begin{proposition}\label{*-invertible}
Let $\ast$ be a star  operation  on  an integral domain  $D$,  and let $ H \in \F(D)$.
 Then  $H$ is $\ast$-invertible if and only if 
$(A:H)^\ast = (A^\ast :H)^\ast = (A H^{-1})^\ast$  for   all $A \in \F(D)$.
\end{proposition}
\begin{proof} Note that, in general,  we have $(A:H)^{\ast }\subseteq (A^{\ast }:H)^\ast = (A^{\ast }:H)$ for all $A, H \in \F(D)$ \cite[page 406, Exercise 1]{G}.

~Assume that $H$ is  $\ast$-invertible,   and
let $x\in (A^{\ast }:H)$.  Therefore $xH\subseteq A^{\ast }$.  Multiplying both sides
by $H^{-1}$ and applying $\ast $,  we get $x\in (A^{\ast }H^{-1})^{\ast
}=(AH^{-1})^{\ast }$.  ~This gives $(A^{\ast
}:H)\subseteq (AH^{-1})^{\ast }$.
Next, let $y \in AH^{-1}$.  Multiplying both
sides by $H$,  we get $yH\subseteq AH^{-1}H\subseteq A$,  and thus $y\in (A:H)$.
So $AH^{-1}\subseteq (A:H)$,  and consequently $(AH^{-1})^{\ast }\subseteq
(A:H)^{\ast }$.   Putting   it all together, we get $(A:H)^{\ast }\subseteq
(A^{\ast }:H)\subseteq (AH^{-1})^{\ast }\subseteq (A:H)^{\ast }$,  which
establishes the equalities.

Conversely, assume that $(A:H)^\ast = (A H^{-1})^\ast$  for all $A \in \F(D)$. In particular, for $A = H$, we have $D \subseteq (H:H)^\ast = (H H^{-1})^\ast \subseteq D$,  and so $H$ is $\ast$-invertible.
\end{proof}

\begin{remark} \rm  Note that  Proposition \ref{*-invertible} can be also deduced from \cite[Corollary 12.1]{hkbook}.
We thank Halter-Koch for pointing out this fact and for informing us that, using the ideal systems  approach on commutative monoids, he has proved a general result on invertibility \cite{hkprivate} that implies the previous Propositions \ref{prop:1.0}, \ref{*CICD}, \ref{*v-cicd},   and   \ref{*-invertible}. 
\end{remark}

 \medskip

We  next give  some examples of $\ast$--CICD's.

\begin{example} \label{ex:1.6} \it  Let $\ast $ be a star operation on an integral domain $D$.  

\bf  Case:  \rm $\ast =v$. \it

  The following properties are    equivalent. 
\begin{enumerate}
\item[(i)] $D$ is a $v$--CICD.
\item[(ii)] $D$ is a $(v, v)$--CICD.
\item[(iii)] $D$ is a CICD.
\item[(iv)] $(AB)^{-1}=(A^{-1}B^{-1})^{v}$   for all  $A,B\in \F(D)$.
 \item[(v)] $D$ is a $v$--multiplication domain. 
\end{enumerate} \rm

The previous statement is an immediate consequence of Propositions   \ref{prop:1.0}, \ref{prop:1.1},  \ref{cor:1.2}, and the fact  that, from the definition of  a  $\ast$--CICD, the notions of  $v$--CICD and  CICD coincide.   Note that the equivalence (iii)$\Leftrightarrow$(v) gives back \cite[Theorem 3.7 ((1)$\Leftrightarrow$(2))]{ET}.  
\smallskip

The case of a star operation of finite character is particularly interesting.  Let  $\ast $ be a star operation on an integral domain $D$. The operation defined by $A^{\ast_{\!{_ f}} } := \bigcup \{ F^\ast \mid F\subseteq A\,,\;  F \in \f(D) \}$  for all $A \in \F(D)$ is a star operation on $D$, called the \it star operation of finite character associated to $\ast$. \rm When $\ast = \ast_{\!{_ f}} $, $\ast$ is called a \it star operation of finite character. \rm As usual, we denote by $t$ the star operation of finite character associated to the $v$--operation, i.e., $t := v _{\!{_ f}}$.   
 We have $\ast_{\!{_ f}} \leq \ast$
for each star operation $\ast$, and    hence, as we have already observed in the   introduction, 
a $\ast_{\!{_ f}}$--CICD is  a   $\ast  $--CICD.  Note that a $\ast_{\!{_ f}}$--CICD  is a special case of  a  $\ast_{\!{_ f}}$--Pr\"ufer domain.   It is obvious from the definitions that the notion of  $\ast_{\!{_ f}}$--Pr\"ufer domain  coincides with that of  \it Pr\"ufer $\ast$--multiplication domain \rm  (for short, P$\ast$MD), i.e., an integral domain such that   $(FF^{-1})^{\ast_{\!{_ f}}} =D$ for   all  $F \in \f(D)$  \cite{hmm}, \cite{FJS}, and \cite{hk}.   

In order to give better interpretations of $\ast_{\!{_ f}}$--CICD's and $(\ast_{\!{_ f}}, v)$--CICD's, we start by 
 recalling that an integral domain $D$ is a \it Dedekind domain \rm  (respectively, { {\it Krull domain}})  if and only if
every $A\in \F(D)$ is  invertible (respectively, $t$-invertible)  (see e.g. \cite[Theorem 37.1]{G} and \cite[Theorem 3.2]{HZ:TV}).    Let $d$ be the identity star operation. Since  $d \leq \ast$ (respectively, $\ast_{\!{_ f}} \leq t$   \cite[Theorem 34.1 (4)]{G})   for  all star operations  $\ast$ ~on $D$, if  $ AA^{-1} =D$  (respectively,  $ (AA^{-1})^{\ast_{\!{_ f}}} =D$), then also $(AA^{-1})^\ast=D$ (respectively, $(AA^{-1})^t=D$).   Therefore a Dedekind domain is a  $ \ast$--CICD for all star operations  $\ast$ on $D$ and a   $\ast_{\!{_ f}}$--CICD is not just a CICD, but more  precisely,   it is a
Krull  domain such  that $A^{\ast_{\!{_ f}}} =A^{t} \ (=A^{v})$ for  all  $A\in \F(D)$ (Proposition \ref{*CICD} ((v)$\Rightarrow$(vi))).

 The previous remarks provide a motivation for the following terminology. Let us call
  a Krull domain such that $A^{\ast_{\!{_ f}}} =A^{t} $ for  all   $A\in \F(D)$   a  \it $\ast $--Dedekind domain. \rm   
Clearly,  a $ \ast_{\!{_ f}}$--CICD coincides with  a  $\ast$--Dedekind  domain   (which is identical by definition   to   a  $ \ast_{\!{_ f}}$--Dedekind domain),  a $v$--Dedekind domain  is just a Krull  domain,   and a  $\ast$--Dedekind  domain   is a particular P$\ast$MD.   
 Next, call a $({\ast_{\!{_ f}}}, v)$--CICD  a     {\it  $(\ast, v)$--Dedekind  domain};    in other words, a \it  $(\ast, v)$--Dedekind  domain \rm   is an integral domain $D$  such that $A^{v}$ is ${\ast_{\!{_ f}}} $-invertible for all $A\in \F(D)$. Obviously the notions of  $(\ast, v)$--Dedekind  domain and $({\ast_{\!{_ f}}}, v)$--Dedekind coincide.  
\end{example}

\begin{example} \label{*=d}   \it  Let $\ast $ be a star operation on an integral domain $D$.  

\bf  Case:  \rm $\ast =d$  (where $d$ is the identity star operation).   

\it The following properties are    equivalent. 
\begin{enumerate}
\item[(i)] $D$ is a      $d$--Dedekind domain   (= $d$--CICD).  
\item[(ii)] $D$ is a Dedekind domain.
\item[(iii)] $\F^v(D)=\F(D)$ and $(AB)^{-1}= A^{-1}B^{-1} $  for  all $A,B\in \F(D)$.

\end{enumerate}
\rm 

As a matter of fact, a Dedekind domain is an integral domain such that every nonzero fractional ideal is invertible (cf. for instance \cite[Theorem 37.1]{G}).  The equivalence of (ii) and (iii) is in  \cite[Corollary 1.3]{Z:G-Dedekind} or \cite[Theorem 2.8]{AK}.   Moreover,   from Proposition \ref{prop:1.1}   we have that   \it the following properties are    equivalent. 
\begin{enumerate}
\item[(j)] $D$ is a   $(d,v)$--Dedekind domain   (= $(d, v)$--CICD).   
\item[(jj)] $D$ is a pseudo-Dedekind domain (i.e., $A^{v}$ is invertible  for   all   $A\in \F(D)$).
\item[(jjj)] $(AB)^{-1}= A^{-1}B^{-1} $  for  all $A,B\in \F(D)$.
\end{enumerate}
 \rm 

Note that      $(d, v)$--Dedekind domains    were
studied under the name of G(enerali\-zed)-Dedekind domains by Zafrullah in 1986 \cite[Theorem 1.1 and Lemma 1.2]{Z:G-Dedekind}
and by D.D. Anderson and Kang  \cite{AK} in 1989 under the name   of   pseudo-Dedekind
domains used above. 
These
domains include locally factorial Krull domains (e.g., UFD's),  rank-one   valuation domains with complete value group,  the ring of entire  functions,   and domains whose groups of
divisibility are complete  lattice-ordered   groups (cf.  \cite[Theorem 1.10, Example 2.1, Theorem 2.6] {Z:G-Dedekind} and \cite[Theorem 2.8 ]{AK}).   If $D$ is a   $(d, v)$--Dedekind  domain,     then  $\F^v(D)$ coincides with the group $\Inv(D)$ of   invertible   ideals  of $D$ (cf. also Corollary \ref{f^v-group} (c));   in the special case where $A^v$ is principal for all $A \in \F(D)$, the set of nonzero fractional $v$-ideals $\F^v(D)$   forms  a group which is isomorphic  to   the group of divisibility of $D$ (Corollary \ref{v-inv-prin}). 

\end{example}

\begin{example} \label{ex:1.7} Let $\ast$ be a star operation  on   an integral domain $D$.  

\bf Case:  \rm $\ast =t$ or $\ast = w$.

\it   The following properties are    equivalent. 
\begin{enumerate}
\item[(i)] $D$ is a   $t$--Dedekind domain   (= $t$--CICD).  
\item[(ii)] $D$ is a    $w$--Dedekind domain  (= $w$--CICD).  
\item[(iii)] $D$ is a Krull domain.
\item[(iv)]  $(AB^{-1} )^{-1}  =(A^{-1}B)^{t}$ for all $A, B \in \F(D)$.  
 \item[(v)]  $(AB )^{-1}  =(A^{-1}B^{-1})^{t}$ and $A^t =A^v$ for all $A, B \in \F(D)$.  
\end{enumerate}
\rm 

\it The following properties are    equivalent. 
\begin{enumerate}
\item[(j)] $D$ is a    $(t, v)$--Dedekind domain  (= $(t, v)$--CICD).  
\item[(jj)] $D$ is a   $(t, w)$--Dedekind domain  (= $(w, v)$--CICD). 
\item[(jjj)]  $D$ is a   pre-Krull domain in the sense of \cite[Proposition 4.1]{Z:ACC}  (i.e., $A^{v}$ is $t$-invertible  for all  $A\in \F(D)$).
\item[(jv)] $(AB)^{-1}=(A^{-1}B^{-1})^{w} =(A^{-1}B^{-1})^{t}$ for all $A, B \in \F(D)$.
\end{enumerate}
 \rm 
 
 The statements (i)--(iii) and the statements (j)--(jv)  are equivalent by  Proposition  \ref{prop:1.1} and from the fact that a $t$-invertible ideal is the same as a $w$-invertible
ideal  \cite[Theorem 2.18]{AC}.    Thus    the  $\ast = w$ case coincides with the  $\ast =t$ case. 

(i)$\Leftrightarrow$(iv)  holds  since (iv) is equivalent to condition (vii) of Proposition \ref{*CICD}   when  $\ast = t$. 

 The fact that (i) implies (v) follows from Propositions \ref{prop:1.1}, \ref{cor:1.2}, and \ref{*CICD} ((v)$\Rightarrow$(vi)).
Conversely, in (v) take  $B:=A^{-1}$;  then  $D\subseteq (AA^{-1} )^{-1} $ ~$ =(A^{-1}A^v)^{t}= (A^{-1}A^t)^{t} = (A^{-1}A)^{t} \subseteq D$, and hence $A$ is $t$--invertible for all $A \in \F(D)$.  

\smallskip

A    $(t, v)$--Dedekind domain   $D$ is a particular Pr\"ufer $v$--multiplication domain (for short, P$v$MD) or, equivalently, a  $t$--Pr\"ufer domain    since $D= (F^v F^{-1})^t = (F^t F^{-1})^t= (FF^{-1})^t $ for all $F\in \f(D)$.    Therefore    $(t, v)$--Dedekind domains    form a class of completely integrally closed  P$v$MD's  that contains the
Krull domains   (and,   \sl a fortiori, \rm   all the $d$--CICD's and the $(d, v)$--CICD's).
Furthermore, we will show (Corollary \ref{f^v-group})   that  for a   $(t, v)$--Dedekind  domain,    the set $\F^v(D)$ of nonzero  fractional   $v$-ideals of $D$ 
 is a complete  lattice-ordered  group under $t$-multiplication.  
\end{example}
\smallskip

The following result is a straightforward
adaptation of  Proposition \ref{*CICD}. 

\begin{proposition} \label{*krull}
\rm  \cite[Theorem 3.9]{AMZ} \it Let $\ast $ be a
star operation on an integral domain  $D$. Then the following  conditions  are  equivalent. 
\begin{enumerate}

\item[(i)]  $(A:B)^{\ast_{\!{_ f}} }=(AB^{-1})^{\ast_{\!{_ f}}}$ for all $A,B\in \F(D)$.

\item[(ii)]  $(A:B^{-1})^{\ast_{\!{_ f}}}=(AB)^{\ast_{\!{_ f}}}$ for all $A,B\in \F(D)$.

\item[(iii)]  $(A^{\ast_{\!{_ f}}}:B)=(AB^{-1})^{\ast_{\!{_ f}}}$ for all $A,B\in \F(D)$.

\item[(iv)]  $(A^{\ast_{\!{_ f}}}:B^{-1})=(AB)^{\ast_{\!{_ f}}}$ for all $A,B\in \F(D)$.

\item[(v)]  $D$ is a   $\ast$--Dedekind domain.  

\item[(vi)] $D$ is a CICD and  $A^{\ast_{\!{_ f}}}=A^t$ for all $A \in \F(D)$. 

\item[(vii)]  $(A^{v}:B^{-1})=(A^{v}B)^{\ast_{\!{_ f}}}$\ for all $A,B\in \F(D)$.
\end{enumerate}
\end{proposition}

The only difference between the above  proposition  and Theorem 3.9 of \cite{AMZ}
is that in \cite{AMZ} what is called here  a   $\ast $--Dedekind domain   is regarded there as a Krull domain, 
which is not correct (i.e., condition (8) of  \cite[Theorem 3.9]{AMZ} is weaker than the other conditions).   For example,  if $\ast =d$,  a  $\ast $--Dedekind domain  is a
Dedekind domain  (Example \ref{*=d}),    which is a very special kind of Krull domain.   This leads to the following  natural    problems.
\begin{problem}
\bf (a) \sl  
 Prove or disprove: There is a finite character star operation $\ast $ that
admits a $\ast $--Dedekind domain $D$ such that {\rm (1)} $D$ is not Dedekind and
{\rm (2)} there is at least one non-Dedekind Krull domain $R$ such that $R$ is not 
$\ast $--Dedekind.

\bf (b) 
 \sl  
Find an example of a   Krull, but not a $\ast$--Dedekind domain, for some $\ast \neq d$,     i.e.,  find a Krull domain  with  a star operation $\ast$ such that   $ d \lneq {\ast_{\!{_ f}}} \lneq t$   (and so $\ast \neq v$).

\end{problem}   
The   difficulty of    problem (a)    lies in the fact   that,   as soon as we  consider the $\ast$-operation on a domain $D$, which is   a  $\ast $--Dedekind  
domain,    this operation    becomes the $t$-operation of the  (Krull)  domain $D$,  as
indicated in (v)$\Leftrightarrow$(vi) of Proposition \ref{*krull}. 
 This problem is important because a positive answer would entail a procedure
for finding finite character operations $\ast $ that admit $\ast $--Dedekind
domains, and the existence of $\ast $--Dedekind domains that answer the
problem would justify deeper study in terms of general $\ast$-operations. 
The negative answer, on the other hand, would give us what we can expect
from a general study.   That is,   we shall know that a $\ast $-Dedekind domain
is either a Krull domain or a Dedekind domain. 

  A positive answer  to    problem (b)
 follows from Example 5.3 of  \cite{FL}, where the authors give the construction of a star operation $\ast$   on a Krull domain  such that $d \lneq \ast = \ast_{\!{_ f}} \lneq t = v$.

\medskip

 Let $\ast$ be a star operation on an integral domain $D$.  Recall that $A\in \F(D)$ is called  {\it  {$\ast$-finite}}  (respectively,  {\it {strictly $\ast$-finite}}) if there exists  an  $F \in  \f(D)$   (respectively,    $F \in \f(D)$ and $F \subseteq A$) such that $A^\ast = F^\ast$ \cite{Z:ACC}. It is well known   that  if   $\ast $ has finite character, then the notions of $\ast$-finite and strictly $\ast$-finite coincide. Moreover, for $A \in \F(D)$,  $A$ is $\ast_f$-invertible if and only if $A$ is $\ast$-invertible and  both  $A$ and $A^{-1}$ are $\ast_f$-finite (for instance   \cite[Lemma 2.3 and Proposition 2.6]{FP},  where this subject was handled in the semistar operation setting).  
\medskip

 The characterizations of   $(t,v)$--Dedekind domains   (or 
 pre-Krull domains)  given in  \cite[Proposition 4.1]{Z:ACC} can be directly translated to the general star  operation  case as
follows. 

\begin{proposition} \label{*-v-dedekind}
Let $\ast $ be a
star operation on an integral domain  $D$. Then the following  conditions  are  equivalent. 
\begin{enumerate}

\item[(i)]  $(AB)^{-1}=(A^{-1}B^{-1})^{\ast_{\!{_ f}}}$   for  all $A, B \in \F(D)$.

\item[(ii)]  $A^{-1}$ is ${\ast_{\!{_ f}}}$-invertible for all $A\in \F(D)$.

\item[(iii)]  $D$ is a  $(\ast, v)$--Dedekind  domain.

\item[(iv)]  $D$ is completely integrally closed and $(AB)^{v}=(A^{v}B^{v})^{\ast_{\!{_ f}}}$
for all $A,B\in \F(D)$.
\end{enumerate}
\end{proposition}
\vskip -0.3cm \begin{proof}
(i)$\Rightarrow$(ii). For all $A \in \F(D)$, clearly we have  $A^vA^{-1} \subseteq D$,   and so $D \subseteq 
(A^vA^{-1} )^{-1}$.  Therefore,  using (i), we have  $D \subseteq 
(A^vA^{-1} )^{-1} = (A^{-1}A^v)^{\ast_{\!{_ f}}}\subseteq D$.  

(ii)$\Rightarrow$(iii). If $A^{-1}$ is  ${\ast_{\!{_ f}}}$-invertible,   then $A^v =(A^{-1})^{-1}$ is also ${\ast_{\!{_ f}}}$-invertible.

 (iii)$\Rightarrow$(i). For all $A, B \in \F(D)$, we have  $A^{-1}B^{-1} \subseteq (AB)^{-1}$,   and so $(A^{-1}B^{-1})^{\ast_{\!{_ f}}} $  $\subseteq ((AB)^{-1})^{\ast_{\!{_ f}}}=(AB)^{-1}$. On the other hand, by assumption, $ D = ((AB)^v (AB)^{-1})^{\ast_{\!{_ f}}} $ and clearly    $((AB)^v (AB)^{-1})^{\ast_{\!{_ f}}} \supseteq (A^vB^v(AB)^{-1})^{\ast_{\!{_ f}}} $;   thus 
$D \supseteq (A^vB^v(AB)^{-1})^{\ast_{\!{_ f}}}$. Multiplying both sides by $A^{-1}B^{-1}$ and applying ${\ast_{\!{_ f}}}$, we have \  $ (A^{-1}B^{-1})^{\ast_{\!{_ f}}} \supseteq $  \, $ (A^vA^{-1}B^vB^{-1}(AB)^{-1})^{\ast_{\!{_ f}}}  = ((AB)^{-1})^{\ast_{\!{_ f}}} = (AB)^{-1}$. We conclude that  $ (AB)^{-1} = (A^{-1}B^{-1})^{\ast_{\!{_ f}}}$.

(iii)$\Rightarrow$(iv). If  $(A^vA^{-1} )^{\ast_{\!{_ f}}} =D$,   then also $(A^vA^{-1} )^t = (A^vA^{-1} )^v =(AA^{-1} )^v  =D$  for  all $A \in \F(D)$. Therefore $D$ is a CICD. For the remainder, since we have already proved (iii)$\Rightarrow$(i), for all $A, B \in \F(D)$  we have $(AB)^v = ((AB)^{-1})^{-1} =  ((A^{-1}B^{-1})^{\ast_{\!{_ f}}})^{-1} = (A^{-1}B^{-1})^{-1} =   ((A^{-1})^{-1}(B^{-1})^{-1})^{\ast_{\!{_ f}}} = (A^vB^v)^{\ast_{\!{_ f}}} $.

(iv)$\Rightarrow$(iii).  Since $D$ is a CICD, for all $A \in \F(D)$, we have $D = (AA^{-1})^v = (A^vA^{-1})^v$. By the equality in (iv), we conclude that $D =  (A^vA^{-1})^v =  (A^v(A^{-1})^v)^{\ast_{\!{_ f}}} $  $ =  (A^vA^{-1})^{\ast_{\!{_ f}}} $. 
 \end{proof}

Obviously, as a $({\ast_{\!{_ f}}}, v)$--CICD is a  $(\ast,
v)$--Dedekind domain,  we can rewrite Proposition \ref{*v-cicd} as a set of  quotient-based 
characterizations of  $(\ast,
v)$--Dedekind domains.   \\

\begin{remark} \label{*-finite}   \bf (a)  \rm
An integral domain $D$ is called \it  $\ast$--Noetherian \rm if $A$ is strictly  $\ast$-finite for all $A\in \F(D)$. It is known that  \it the following  conditions  are equivalent.  
\begin{enumerate}
\item[(i)]  $D$ is  $\ast$--Noetherian.
\item[(ii)]  $D$ is  $\astf$--Noetherian.  
\item[(iii)]  $D$ satisfies the ascending chain condition on $\ast$-ideals.
\end{enumerate} \rm
For the proof of the previous statement and more details on this  subject,  cf. \cite[Lemma 3.3 and Proposition 3.5]{EFP}.  
Note also that, if $D$ is  $\ast$--Noetherian,  then each ideal of $D$ is $\ast$-finite, but the converse is false in general   (cf.    \cite[page 29]{DDA-DFA}  and  \cite[Example 18]{GJS}).  

 From the previous   considerations,   we easily deduce that a  $\ast$--Dedekind domain is  a $\ast$--Noetherian domain.   Moreover, as observed above, a $\ast$--Dedekind domain is a P$\ast$MD.  Note that the converse is also true, i.e., the $\ast$--Dedekind domains coincide with the $\ast$--Noetherian P$\ast$MD's.  As a matter of fact,     
if $D$ is $\ast$--Noetherian,   then  for  all   $A \in \F(D)$, there  exists  an   $F \in \f(D)$   with $F \subseteq A$     such that  $F^\ast = A^\ast = A^{\astf}$. 
 Hence $F^v = (F^\ast)^v = (A^{\ast})^v =A^v$; thus $F^{-1} = A^{-1}$.   Therefore 
$(AA^{-1})^{\astf} = (A^{\astf} A^{-1})^{\astf} = (F^{\astf} F^{-1})^{\astf} = (F F^{-1})^{\astf} $.  If we assume that $D$ is also a P$\ast$MD,  then  $(FF^{-1})^{\astf} = D$;   hence  $(AA^{-1})^{\astf}=D$, i.e., $D$ is $\ast$--Dedekind.   

 From the previous observations we can conclude  that  the  notion of $\ast$--Dedekind domain, given here, coincides in the star operation case with the notion considered for semistar operations  in \cite[Proposition 4.1]{EFP}.

We can summarize some of the previous considerations by saying that \it  the following notions coincide.
\begin{enumerate}
\item[(i)]  $\ast$--Dedekind \  (= $\astf$--Dedekind = $\astf$--CICD).
\item[(ii)] $\astf$--multiplication domain.
\item[(iii)] $\ast$--Noetherian and Pr\"ufer $\ast$--multiplication domain.
\item[(iv)] $\astf$--Noetherian and $\astf$--Pr\"ufer domain.
\end{enumerate}

\bf (b) \rm  Note that Noetherian ideal systems are  investigated in  \cite[Chapter 3]{hkbook}. In  particular,  the equivalent statements given in (a) are also proved in \cite[Theorem 3.5]{hkbook} in the more general setting of ideal systems on monoids.

Given an ideal system $r$, $r$--Dedekind monoids are introduced  and studied in \cite[Chapter 23, \S3]{hkbook}. 
 However, the  notion of $\ast$--Dedekind domain  coincides   with that in \cite{hkbook} in case $\ast= \astf$ or in  the   case of Krull domains, but they are different in general (e.g., $v$--Dedekind domains are precisely   Krull domains,   but  $v$--Dedekind monoids are just completely integrally closed monoids).

We  also  note that, using the ideal systems  approach on commutative monoids, Halter-Koch \cite{hkprivate} has obtained  a general version of  Proposition \ref{*-v-dedekind}. 

\end{remark}

 
\section{Star Pr\"ufer domains }

Recall that an integral domain  $D$ is a \it $v$--domain \rm if each $F\in \f(D)$ is $v$-invertible.  We have already introduced a 
direct generalization of this definition    when   $\ast$ is a star operation on $D$   by  saying  that $D$ is a \it   $\ast$--Pr\"ufer  domain   \rm     if every $F \in \f(D)$ is $\ast $-invertible. 
Since  a $\ast $-invertible ideal  is always $v$-invertible, we observe that a   $\ast$--Pr\"ufer domain    is always a $v$--domain.    Note that    $\ast$--Pr\"ufer  domains  were recently introduced  in  the  case of semistar operations $\star$ under the name of $\star$--domains \cite[Section 2]{FP2}.

If $F \in \f(D)$ is  $\ast $-invertible,   then $
F^{\ast }=F^{v}.$ Since, in a   $\ast$--Pr\"ufer domain,   this holds for  all  
$F\in \f(D)$, we conclude that  ${\ast_{\!{_ f}}}=t$ in   a   $\ast$--Pr\"ufer domain.

 Next, we can consider a weaker notion:  call $D$ a \it  $(\ast, v)$--Pr\"ufer domain   \rm if $F^{v}$ is $\ast $-invertible for all $F\in \f(D)$.  It is easy to see that a   $(\ast, v)$--Pr\"ufer domain    is also a $v$--domain and that a   $\ast $--Pr\"ufer domain    is a  
 $(\ast, v)$--Pr\"ufer domain.        Clearly, if  $\ast_1, \ast_2$ are two star operations on $D$ and $\ast_1 \leq  \ast_2$, then a   $\ast_1 $--Pr\"ufer domain (respectively,  $(\ast_1, v) $--Pr\"ufer domain)   $D$ is a $\ast_2 $--Pr\"ufer  domain  (respectively, $(\ast_2, v) $--Pr\"ufer  domain).    

 We have already observed that, from the definitions, it follows immediately  that  the notions of   ${\ast_{\!{_ f}}}$--Pr\"ufer  domain    and  P$\ast$MD (or  P${\ast_{\!{_ f}}}$MD) coincide. Therefore a P$\ast$MD is a   $\ast $--Pr\"ufer domain,  but the converse is not true since there are $v$-domains (=  $v$--Pr\"ufer domains) that are not P$v$MD's \cite{HO}. \  
  Also note   that for an ideal system $r$ on a monoid, the notion of $r$--Prufer
monoid, for a general $r$, was introduced in \cite[Chapter 17]{hkbook}. However, most
of the results on $r$--Prufer monoids in \cite{hkbook} were proved for $r$-finitary.
Now $r$--Prufer monoids for finitary $r$ coincide with $\ast $--Pr\"ufer domains
only in case $\ast ={\ast_{\!{_ f}}}$. That is, the $r$--Pr\"ufer monoids studied in
\cite{hkbook} were simply P$\ast $MD's in ring-theoretic terms.

\begin{example} \label{ex:2.1} Let $\ast$ be a star operation defined on an integral domain $D$.  

\bf Case:  \rm $\ast =d$.

Clearly, from the definition, \it  the following notions  coincide.  
\begin{enumerate}
\item[(i)] $D$ is a   $d$--Pr\"ufer domain.  
\item[(ii)] $D$ is a Pr\"ufer domain. 
\item[(iii)] Each $F \in \f(D)$ is invertible.
\end{enumerate} \rm

\it  The following notions  coincide.  
\begin{enumerate}
\item[(i)] $D$ is a   $(d, v)$--Pr\"ufer  domain.   
\item[(ii)] $D$   is a   generalized  GCD (for short, GGCD) domain (i.e., the intersection of two invertible ideals is invertible).
\item[(iii)]  $F^v$ is invertible for all   $F \in \f(D)$.
\end{enumerate} \rm 

   Generalized  GCD 
domains were introduced in   \cite{AA},    where the previous equi\-va\-lence was also proven  \cite[Theorem 1]{AA}.

Note that,  while a Pr\"{u}fer domain is a GGCD domain,  there
are examples of GGCD domains that are not Pr\"{u}fer  \cite[Theorem 2 (2)]{AA}. 
So, while a 
 $\ast $--Pr\"ufer  domain is a $(\ast, v)$--Pr\"ufer domain, a $(\ast, v)$--Pr\"ufer domain may not be a  $\ast $--Pr\"ufer domain.

\bf Case:  \rm $\ast =t$ or $\ast =w$.

\it   The following notions  coincide. 
\begin{enumerate}
\item[(i)] $D$ is a   $t$--Pr\"ufer domain.  
\item[(ii)] $D$ is a  $(t, v)$--Pr\"ufer domain.  
\item[(iii)] $D$ is a   $w$--Pr\"ufer domain.   
\item[(iv)] $D$ is a   $(w, v)$--Pr\"ufer domain.  
\item[(v)] $D$ is  Pr\"ufer $v$--multiplication domain.  
\end{enumerate}
\rm 
Since the maximal $t$-ideals coincide with the maximal $w$-ideals, the notions  $w$-invertible and $t$-invertible coincide (cf.   \cite[Theorem 2.18]{AC}  and \cite[Section 5]{WM}), thus   (iii)$\Leftrightarrow$(i)$\Rightarrow$(ii)$\Leftrightarrow$(iv)   and, by definition, (v)  coincides  with (i).  Finally, (ii)$\Rightarrow$(v), since $D = (F^vF^{-1})^t = (F^tF^{-1})^t = (FF^{-1})^t $ for all $F \in \f(D)$.

 \bf Case:  \rm $\ast =v$.

From the definitions, we  immediately have  that  \it the following notions   coincide.  
\begin{enumerate}
\item[(i)]  $D$ is a    $v$--Pr\"ufer domain.  
\item[(ii)]  $D$ is a   $(v,v)$--Pr\"ufer  domain.  
\item[(iii)]  $D$ is a $v$--domain.
\end{enumerate} \rm 
\end{example}

\bigskip   Once we know that a Pr\"{u}fer domain is just a special case of a    $\ast $--Pr\"ufer  
domain  (and, in particular, of  a $(\ast, v) $--Pr\"ufer  
domain) we would like to see what ideal-theoretic   characterizations  of Pr\"{u}fer
domains can be translated to the framework of $\ast $--Pr\"ufer and $(\ast, v) $--Pr\"ufer domains.   
Here we point out some.

\begin{theorem} \label{*prufer} Let $\ast$ be a star operation defined on an integral domain $D$.  Then the   following properties are  equivalent.  

\begin{enumerate}
\item[(i)]  $D$ is a  $\ast$--Pr\"ufer domain. 
\item[(ii)]  Every (nonzero) two generated
ideal of $D$ is $\ast $-invertible.
\item[(iii$_{\f}$)]  $((F\cap
G)(F+G))^{\ast }=(FG)^{\ast }$ for all $F,G \in \f(D)$.

\item[(iii$_{\F}$)]   $((A\cap
B)(A+B))^{\ast }=(AB)^{\ast }$ for all $A,B\in \F(D)$.

\item[(iv$_{\f}$)] $(F(G^{\ast }\cap
H^{\ast }))^{\ast }=(FG)^{\ast }\cap (FH)^{\ast }$ for all $F, G, H\in \f(D)$.
\item[(iv$_{\!\f\!\F}$)] $(F(A^{\ast }\cap
B^{\ast }))^{\ast }=(FA)^{\ast }\cap (FB)^{\ast }$ for all $F\in \f(D)$ and $
A,B\in \F(D)$.

 \item[(v)] If  $A, B \in \F(D)$ are $\ast$-invertible, then $A \cap B$ and $A+B$ are $\ast$-invertible.  

 \item[(vi)] If  $A, B \in \F(D)$ are $\ast$-invertible, then $A+B$ is $\ast$-invertible.  
\end{enumerate}

\end{theorem} 

\begin{proof}
The proof of (i)$\Leftrightarrow$(ii)   follows   from   the reduction argument used in  \cite[Lemma 2.6]{MNZ} for showing that an integral domain is a $v$--domain if and only if   every (nonzero) two generated
ideal is $v$-invertible.  
 For the sake of completeness, we give some details of the proof of (ii)$\Rightarrow$(i).  Let $F \in \f(D)$;   we want to show that $F$ is $\ast$-invertible.  We use induction on the number of generators of $F$. Let $F :=  (x_1 , x_2,
..., x_n, x_{n+1})D$ with $n \geq 2$ and set $I :=x_1D$,   $J :=  (x_2, x_3,  ..., x_n)D$,   and $H :=
x_{n+1}D$. Then $F(IH + IJ +JH) = (J + H)(H + I)(I + J)$. Note that each of
the factors on the right is generated by \    $ k  \leq  n$ \  elements,
and so is $\ast$-invertible by the induction hypothesis. This forces the factors on the left (and hence, in particular, $F$) to
be $\ast$-invertible. \  (Note
that this method of proof is essentially that used originally  by H. Pr\"{u}fer in \cite[page 7]{P}
to show that $D$ is a Pr\"{u}fer domain if and only if every (nonzero) two generated
ideal of $D$ is invertible.)

(i)$\Rightarrow$(iii$_{\f}$).   $((F\cap
G)(F+G))^{\ast }\subseteq (FG)^{\ast }$ holds for all    $F,G\in \f(D)$.

For the
reverse  containment,    let $x\in FG$. Then $xG^{-1}\subseteq FGG^{-1}\subseteq F$ and $
xF^{-1}\subseteq F^{-1}FG\subseteq G$.  This gives $x(F^{-1}\cap
G^{-1})\subseteq F\cap G$.   But  $F^{-1}\cap G^{-1}=(F+G)^{-1}$.  So we have $
x(F+G)^{-1}\subseteq F\cap G$.  Multiplying both sides by  $F+G$  and  applying  $\ast$,   we get $x\in ((F\cap G)(F+G))^{\ast }$. This gives $
(FG)^{\ast }\subseteq ((F\cap G)(F+G))^{\ast }$.

(iii$_{\f}$)$\Rightarrow$(ii) is
obvious because   $((F\cap G)(F+G))^{\ast }=(FG)^{\ast }$  for  all $F, G\in \f(D)$
implies that in particular $((xD\cap yD)(xD +yD))^{\ast }=xyD$ for all  nonzero $
x, y\in D$,  which forces every (nonzero)  two generated  ideal of 
$D$ to be $\ast$-invertible.

(iii$_{\f}$)$\Rightarrow$(iii$_{\F}$). Obviously,   $((A\cap B)(A+B))^{\ast }\subseteq (AB)^{\ast }$ holds  for all $A,B\in
\F(D)$.    For the reverse 
containment,   it is enough to show that $AB\subseteq ((A\cap B)(A+B))^{\ast }$.
For  this,   let $x\in AB$.  Then  $x\in FG$,   where $F$ and $G$ are finitely
generated  with  $F\subseteq A$ and $G\subseteq B$.  But then,  by (iii$_{\f}$),
$x\in (FG)^{\ast }=((F\cap G)(F+G))^{\ast }\subseteq
((A\cap B)(A+B))^{\ast }$. Thus $AB\subseteq ((A\cap B)(A+B))^{\ast }$.

(iii$_{\F}$)$\Rightarrow$(iii$_{\f}$)  and (iv$_{\!\f\!\F}$)$\Rightarrow$(iv$_{\f}$) are trivial.

(i)$\Rightarrow$(iv$_{\!\f\!\F}$). Obviously $(F(A^{\ast }\cap B^{\ast }))^{\ast }\subseteq
(FA)^{\ast }\cap (FB)^{\ast }$.   For the reverse  containment,   note that $F$ is 
$\ast $-invertible. Now consider $(F^{-1}((FA)^{\ast } \cap (FB)^{\ast
}))^{\ast }\subseteq (F^{-1}(FA)^{\ast } \cap F^{-1}(FB)^{\ast })^{\ast
}\subseteq (F^{-1}(FA))^{\ast } \cap (F^{-1}(FB))^{\ast }$ $= A^{\ast } \cap
B^{\ast }$.   So the inclusion $(F^{-1}((FA)^{\ast }$ $\cap (FB)^{\ast }))^{\ast } \subseteq
A^{\ast }\cap B^{\ast }$   gives, on multiplying by $F$ and applying $
\ast $ on both sides,  the reverse containment.

(iv$_{\f}$)$\Rightarrow$(ii).
Let $F:=(a,b)$, 
  $
G:=(\frac{1}{a})$ and $H:=(\frac{1}{b})$, where $a$ and $b$ are two nonzero elements of $D$. Then, by assumption,  we have $((a, b)((\frac{1}{a})\cap (\frac{1}{b}))^{\ast }=((\frac{1}{a})(a,b))^{\ast }\cap ((\frac{1}{b})(a,b))^{\ast }$. On the other hand, it is easy to see that $(a,b)(a,b)^{-1} = (a, b)((\frac{1}{a})\cap (\frac{1}{b}))$.  Therefore   $D \supseteq  ((a,b)(a,b)^{-1} )^\ast = ((\frac{1}{a})(a,b)))^{\ast }\cap ((\frac{1}{b})(a,b))^{\ast } \supseteq D$, 
and so we conclude that $((a,b)(a,b)^{-1})^{\ast }=D$.

(v)$\Rightarrow$(vi)$\Rightarrow$(ii) are obvious (for the last implication note that a nonzero principal ideal is ($\ast$-)invertible).

(iii$_{\F}$)$\Rightarrow$(v).  Since $A, B \in \F(D)$ are  $\ast$-invertible if and only if $AB$ is $\ast$-invertible, the  conclusion follows from the equality $((A\cap
B)(A+B))^{\ast }=(AB)^{\ast }$.

\end{proof}
%

\begin{remark}  \rm
Let $D$ be a $\ast$--Pr\"ufer
domain. If we assume that $\ast $ has finite character (hence, $D$ is a P$\ast$MD), then  (as we observed above)  $A\in \F(D)$ is $\ast $-invertible if and only if $A$ is  (strictly)  $\ast$-finite.  In this case, (vi) of Theorem
\ref{*prufer}  reduces to  ``the sum of two $\ast$-finite ideals is $\ast$-finite".
  (This is a  trivial statement since $(F^\ast+G^{\ast})^\ast = (F+G)^\ast$ for all $F,G \in \f(D)$,  \cite[Proposition 32.2]{G}.)  
 
However, if $\ast $\ does not have finite character, a $\ast $-invertible
ideal need not be strictly $\ast$-finite. In fact   (for $D$  a  $\ast $--Pr\"ufer domain),  each $\ast 
$-invertible ideal is strictly $\ast$-finite precisely when $D$ is  a   $\ast _{f}$--Pr\"ufer domain (= P$\ast$MD).
 (If $F\in \f(D)$ is $\ast$-invertible,   then $F^{-1}$ is  $\ast$-invertible,    and so there  exists  $G \in \f(D)$ such that  $G \subseteq F^{-1}$ and $F^{-1} =(F^{-1})^\ast = G^\ast$. Therefore   $D = (FF^{-1})^\ast =  (F(F^{-1})^\ast)^\ast = (FG^\ast)^\ast =  (FG)^\ast = (FG)^{\ast_f} = (FF^{-1})^{\astf} \subseteq D$; hence $(FF^{-1})^{\astf}= D$.) 
\end{remark}

It is natural to ask whether it is possible  to remove from   statement  (iv$_{\!\f\!\F}$) of the previous theorem  the
condition that $F$ is finitely generated.  
We do not have a complete answer to this question,
however  we have an interesting alternative described in the following proposition. 
Recall that a star operation $\ast 
$ is called \it stable \rm (or, \it distributes  over finite intersections\rm ) if $(A\cap B)^{\ast }=A^{\ast }\cap B^{\ast }$ for all $
A,B\in \F(D)$ (cf. \cite[page 174]{FH} and \cite{AC-2005}).  A star   operation   induced by a defining family of quotient
rings of $D$ is stable (cf. for instance \cite[Proposition 2.2]{AC-2005}).

\begin{proposition}
Let $\ast$ be a star operation defined on an integral domain $D$.  Then the   following properties are  equivalent.  

\begin{enumerate}
\item[($\overline{\mbox{i}}$)]  $D$ is a   $\ast$--Pr\"ufer domain     and $\ast$ is a stable star operation on $D$.

\item[($\overline{\mbox{iv}}$)] $\left(C(A \cap B)\right)^{\ast} = \left(C(A^{\ast }\cap B^{\ast })\right)^{\ast }=(CA)^{\ast }\cap (CB)^{\ast }$ for all $A,B, C \in \F(D)$.
\end{enumerate}
\end{proposition}
\begin{proof}
($\overline{\mbox{i}}$)$\Rightarrow$($\overline{\mbox{iv}}$).  As seen above, in general we have $(C(A^{\ast }\cap B^{\ast }))^{\ast }\subseteq
(CA)^{\ast }\cap (CB)^{\ast }$.  Moreover, since $\ast$ is stable,   
$\left(C(A \cap B)\right)^{\ast} = \left(C(A \cap B)^{\ast }\right)^{\ast} = (C(A^{\ast }\cap B^{\ast }))^{\ast }$. \ 
 For the reverse  
containment,  it is sufficient to show that $CA\cap CB\subseteq (C(A\cap B))^{\ast }$.   For this,  
let $x\in CA\cap CB$.  Then, in particular,  there is  an  $F\in \f(D)$ such that $x\in FA\cap FB$ and $F \subseteq C$.
So $x\in (FA\cap FB)^{\ast }=(FA)^{\ast }\cap (FB)^{\ast }
 =(F(A^{\ast
}\cap B^{\ast }))^{\ast }$  (the last equality holds by Theorem \ref{*prufer} ((i)$\Rightarrow$(iv$_{\!\f\!\F}$))).  Again,    by the stability of  $
\ast $,   we have $(F(A^{\ast }\cap B^{\ast }))^{\ast }=(F(A\cap B))^{\ast }$.
Thus we conclude that $x\in (FA\cap FB)^{\ast }=(F(A\cap B))^{\ast }\subseteq
(C(A\cap B))^{\ast }$.

($\overline{\mbox{iv}}$)$\Rightarrow$($\overline{\mbox{i}}$). From $(C(A\cap B))^{\ast }=(CA)^{\ast }\cap (CB)^{\ast }$
for all  $A,B,C\in \F(D)$,   by setting $C=D$,  we deduce that $\ast $ is stable.  Moreover, as in the proof of Theorem \ref{*prufer} ((iv$_{\f}$)$\Rightarrow$(ii)),
taking $C:=(a,b)$, 
   $A:=(\frac{1}{a})$,   and $B:=(\frac{1}{b})$, where $a$ and $b$ are two nonzero elements of $D$, we obtain that $((a,b)(a,b)^{-1})^{\ast }=D$.
\end{proof}

As we have already observed, the notions of    ${\ast_{\!{_ f}}}$--Pr\"ufer  domain    and  P$\ast$MD (or  P${\ast_{\!{_ f}}}$MD) coincide;  moreover,  for a
P$\ast$MD  the operation ${\ast_{\!{_ f}}}$ is stable \cite[Theorem 3.1]{FJS}.   Furthermore, a $\ast$--CICD is a particular  $\ast$--Pr\"ufer  domain,     and for a $\ast$--CICD the operation $\ast$ is stable \cite[Theorem 2.8]{AC-2005}.  Therefore,   in order to find an example of a  ~$\ast$--Pr\"ufer domain  for which  $\ast $ is not  stable,   we have to consider the case of a  $\ast$--Pr\"ufer domain,   not   a  $\ast$--CICD, with  a   star operation $\ast$  that is not of finite character.    In case $\ast =v$,   it is known that $D$ is a $v$--domain if and only if $D$ is integrally closed  and   $v$ distibutes over finite intersections of finitely generated ideals  (cf. \cite[Theorem 1(2) and Theorem 2]{MO}   or \cite[Theorem 2.8]{AC-2005}).  Very  recently,   Mimouni has given an explicit example of a two-dimensional Pr\"ufer domain (hence a   $v$--domain,  but not  a  ($v$--)CICD)
with $v$ not stable \cite[Example 3.1]{Mi}.

\begin{proposition} \label{*prufer-necessary}  Let $\ast$ be a star operation defined on an integral domain $D$.   If   $D$ 
is a  $\ast$--Pr\"ufer domain,      then:
\begin{enumerate}
\item[(1)] $(FG)^{-1}=(F^{-1}G^{-1})^{\ast }$ for all $F, G\in \f(D)$.

\item[(2)]  If $A, B \in \F(D)$ are $\ast $-invertible,  then $A\cap B$ is $\ast $-invertible.  

\item[(3)] $\ast$ is an a.b. star operation (i.e.,  $(FA)^\ast \subseteq
(FB)^{\ast }$ implies that $A^{\ast }\subseteq B^{\ast } $  for   all $F\in
\f(D) $ and $A, B \in \F(D)$). In particular, $D$ is an  integrally closed domain \cite[Corollary 32.8]{G}.

\end{enumerate}
\end{proposition}

\begin{proof} (1).    If $F,G \in \f(D)$, then in particular $G^{-1}$ is $\ast$-invertible. By Proposition \ref{*-invertible},  $(F^{-1}G^{-1})^{\ast } = (F^{-1} : G)^\ast$. The conclusion follows from the fact that $(F^{-1} : G)^\ast = ((D:F) :G)^\ast = (D: FG)^\ast = ((FG)^{-1})^\ast = (FG)^{-1}$.

 That (2) holds was already observed in
Theorem \ref{*prufer} ((i)$\Rightarrow$(v)).

(3). The proof is  straightforward;   multiply both sides of the relation $(FA)^\ast \subseteq
(FB)^{\ast }$ by $F^{-1}$ and apply the $\ast $-operation.
\end{proof}

\begin{remark} \label{rk:2.5}
\rm  \bf (a) \rm To see that  statement (1) of  Proposition \ref{*prufer-necessary}   does not   characterize $\ast$--Pr\"ufer domains,    recall that an integral
domain $D$ is a \it  pre-Schreier domain \rm  if for all nonzero $x,y,z\in D$,  
 $x|yz$   implies that   $x=rs$,  where  $r|y$ and $s|z$.  Pre-Schreier
domains are a generalization of GCD domains (cf. \cite{Cohn} and \cite[Theorem 1]{Dribin}).  It   is  well known  that a pre-Schreier domain   satisfies    (each of)  the following equivalent  properties.  
\begin{enumerate}
\item[\bf ($\boldsymbol  \alpha$)\rm ]
 For all $0 \neq a_i, b_j \in D$, with $1\leq i \leq n, \  1\leq j \leq m$:
$$
\left(\cap_{1\leq i \leq n} (a_{i})\right)
\left(\cap_{1\leq j \leq m} (b_{j})\right) =
\cap_{1\leq i \leq n, 1\leq j \leq m}   (a_{i}b_{j}).
$$

\item[\bf ($\boldsymbol \beta$)\rm ] $(FG)^{-1}=F^{-1}G^{-1}$ for all $F,G\in \f(D)$.
\end{enumerate}
(Cf. \cite[Corollary 1.7]{Z:pre-Schreier}  and
 \cite[Proposition 1.6]{Z:G-Dedekind};  note that an integral domain satisfying these two equivalent conditions was called  a $\ast$-domain  by Zafrullah in \cite{Z:pre-Schreier}).   Now there do exist
pre-Schreier domains that are not integrally closed \cite[page
1918]{Z:pre-Schreier}.  Combining this piece of information 
with the   fact that  a    $\ast$--Pr\"ufer domain    is integrally closed (Proposition  \ref{*prufer-necessary} (3)), 
we easily conclude that $(FG)^{-1}=(F^{-1}G^{-1})^{\ast }$ for all $F,G\in \f(D)$ does not imply that $D$ is  a  $\ast$--Pr\"ufer domain.

\bf (b) \rm  Statement (2) of Proposition \ref{*prufer-necessary}   does not characterize $\ast$--Pr\"ufer domains.  For example,
in a  generalized  GCD domain (Example \ref{ex:2.1}),    we have that \ $A\cap B$ \  is   
($d$-)invertible for all ($d$-)invertible  $A, B \in \F(D)$.     However, a  generalized  GCD domain  may not be  a   $d$--Pr\"ufer domain    (= Pr\"{u}fer domain, Example \ref{ex:2.1})  \cite[Theorem 2 (2)]{AA}.   \ (It may also be
noted that for a mere pair of ideals   $F,G\in \f(D)$,   the fact that $F\cap G$ \ is $
\ast $-invertible does not mean that $F$ and $G$ are both 
$\ast$-invertible.  For instance, let $k$ be a field and $X, Y$ two indeterminates over  $k$.  Take  $ D:=k[X,Y]$, $F:= (X, Y)$ and $G:=(X)$, then  $(X,Y)\cap (X)$ is principal,  but $(X,Y)$
is not invertible.)

\bf (c)  \rm  Statement (3) of Proposition \ref{*prufer-necessary}  does not characterize $\ast$--Pr\"ufer domains    (cf. \cite[Example 1  (2),  page 150]{FP2}).  
 It is easy to show that an integral domain  $D$ is 
integrally closed   if and only if  there is an a.b. star
operation $\ast $  defined on $D$.   (The proof depends upon the fact that if $D$ is integrally
 closed,   then $D$ is expressible as an intersection of valuation overrings of $
D$ \cite[Theorems 19.8, 32.5 and Corollary 32.8]{G}.) 
 Therefore,   an integrally closed domain may not be  a   $\ast $--Pr\"ufer domain   for
any star operation $\ast $  \ since   there are
integrally closed domains that are not $v$--domains  (e.g.,
\cite[page 429, Exercise 2]{G}).

\end{remark}

Bearing in mind  Proposition \ref{*prufer-necessary}  (1),   we  next  give   more precise relations among the notions coming into play.

\begin{proposition} \label{pr:2.6}  Let $\ast$ be a star operation defined on an integral domain  $D$,   and consider the following statements:
\begin{enumerate}
\item[(a)] $D$ 
is a   $\ast$--Pr\"ufer domain.    

 \item[(b)]  $(AF)^{-1}=(A^{-1}F^{-1})^{\ast }$ for all $
A \in \F(D)$ and $F\in \f(D)$.

\item[(c)]  $D$ is a   $(\ast, v)$--Pr\"ufer domain.   
\end{enumerate}
Then \rm  (a)$\Rightarrow$(b)$\Leftrightarrow$(c). 
\end{proposition}

\begin{proof} 
 (a)$\Rightarrow$(b). Let $
A \in \F(D)$,   $F\in \f(D)$,   and let  $x\in A^{-1}F^{-1}$. Then  
$xF\subseteq A^{-1}F^{-1}F \subseteq A^{-1}$,   and so $xAF\subseteq
A^{-1}A\subseteq D$.  Therefore $A^{-1}F^{-1}\subseteq (AF)^{-1}$.  On the other  hand,   
$(AF)^{-1}$ is a $v$-ideal  (and so, in particular, a   $\ast$-ideal);  thus we have $(A^{-1}F^{-1})^{\ast }\subseteq ((AF)^{-1})^\ast =
(AF)^{-1}$.   Conversely,  let $y \in (AF)^{-1}$. Then ~$yAF\subseteq (AF)^{-1}AF\subseteq D$, so 
$yF\subseteq A^{-1}$.  Multiplying both sides by $F^{-1}$, applying $\ast ,$
and noting that $F\in \f(D),$ we get $y \in y (FF^{-1})^\ast \subseteq (A^{-1}F^{-1})^\ast$.

(b)$\Rightarrow$(c).    Let  $F \in \f(D)$,   and set $A:=F^{-1}$. Then   $D \subseteq (F^{-1}F)^{-1} = (AF)^{-1}=(A^{-1}F^{-1})^{\ast } = (F^vF^{-1})^{\ast } \subseteq D$,  and so $D$ is a   $(\ast, v)$--Pr\"ufer domain.

  (c)$\Rightarrow$(b). Let $ x \in (AF)^{-1}$. We have  $xAF \subseteq D$;    so $xF
\subseteq A^{-1}$,   and hence $xF^v = (xF)^v \subseteq (A^{-1})^v = A^{-1}$. Therefore $xF^vF^{-1} \subseteq A^{-1}F^{-1}$. Since $ (F^vF^{-1})^\ast = D$, we conclude that $(AF)^{-1} \subseteq  (A^{-1}F^{-1})^\ast$. 
The reverse containment is straightforward since $((AF)^{-1})^\ast = (AF)^{-1}$.  
\end{proof}

 Note that we already observed that there are  $(d,v)$--Pr\"ufer domains    that are not   $d$--Pr\"ufer domains  (= Pr\"ufer domains,  Example \ref{ex:2.1}),   and so   (c) of Proposition  \ref{pr:2.6}   does  not imply (a).   We will see later (Theorem \ref{*-inv} (c)) that (c) and (a) are equivalent under   an additional  condition. 
\medskip 

For   $\ast$--Pr\"ufer  domains,     we have the following set of  quotient-based  characterizations.

\begin{theorem} \label{*prufer-2} Let $\ast$ be a star operation defined on an integral domain $D$.   Then the following conditions  are   equivalent.  
\begin{enumerate}
\item[(i)] $D$ is  a  $\ast$--Pr\"ufer domain.   
\item[(ii)] For  all  $A\in \F(D) $ and $F\in \f(D)$, $A\subseteq F^{\ast }$ implies $A^{\ast
}=(BF)^{\ast }$ for some $B\in \F(D)$.

\item[(iii)] $(A:F)^{\ast }=(A^{\ast }:F) = (AF^{-1})^{\ast }$  for all  $A\in \F(D) $ and $F\in \f(D)$.

\item[(iv)]  $(A:F^{-1})^{\ast }= (A^{\ast }:F^{-1}) = (AF)^{\ast }$ for all  $A\in \F(D) $ and $F\in \f(D)$.

\item[(v)]  $(F:A)^{\ast }=(F^{\ast }:A)=(FA^{-1})^{\ast }$ for all  $A\in \F(D) $ and $F\in \f(D)$.

\item[(vi)]  $(F:A)^{v}= (F^{v}:A)=(FA^{-1})^{\ast }$ for all  $A\in \F(D) $ and $F\in \f(D)$.

\item[(vii)]  $(F^{v}:A^{-1})=(FA^{v})^{\ast }$ for all  $A\in \F(D) $ and $F\in \f(D)$.

\item[(viii)] $((A+B):F)^{\ast}=((A:F)+(B:F))^{\ast }$  for all $A,B\in \F(D)$  and   $F\in \f(D)$.

 \item[(ix)]  $(A :(F\cap G))^\ast = ((A:F) + (A:G))^\ast$  for   all  $A\in \F(D) $ and $F, G\in \f^\ast(D):= \{ H \in \f(D) \mid H = H^\ast \}$.  

 \item[(x)]   $(((a):_{D}(b))+ ((b):_{D}(a)))^{\ast }=D$   for all  nonzero  $
a,b\in D$.
\end{enumerate}
\end{theorem}

\vskip 20pt
 \begin{proof}  (i)$\Rightarrow$(ii). Set $B:=(AF^{-1})^\ast$. Then clearly $(BF)^\ast = ((AF^{-1})^\ast F)^\ast =  (A (FF^{-1})^\ast)^\ast = A^\ast$.

(ii)$\Rightarrow$(i). We
show that every $F\in \f(D)$ is $\ast $-invertible.   For this,  let $0\neq x \in F^{\ast }$,   and set $A := (x)$.  Then by  assumption,  there is  a  $B\in \F(D)$
such that $(x) = (x)^\ast = A^\ast =(BF)^{\ast }$,  and so $D = ((x^{-1}B)F)^\ast$,   which is equivalent to $F$ being $\ast $-invertible.

(i)$\Rightarrow$(iii)  follows from Proposition \ref{*-invertible}, since in the present situation $F \in \f(D)$ is $\ast$-invertible.

(i)$\Rightarrow$(iv). This implication can be proven in a similar fashion as (i)$\Rightarrow$(iii)
using the fact that if $F$ is  $\ast $-invertible, then  so   is $F^{-1}$.

(i)$\Rightarrow$(v).  Clearly $(F:A)^{\ast }\subseteq (F^{\ast
}:A)$.
 Let $x\in
(F^{\ast }:A)$. Then  $xA\subseteq F^{\ast }$;   so $xAF^{-1}\subseteq F^{\ast
}F^{-1} \subseteq (F^{\ast
}F^{-1})^\ast = D$,  which gives $xF^{-1}\subseteq
A^{-1}$.  Now multiplying both sides by $F$ and applying $\ast $,  we get $x\in
(FA^{-1})^{\ast }$.
  Next,  to show that $(FA^{-1})^{\ast }\subseteq
(F:A)^{\ast },$ let $y\in FA^{-1}$.  Then  $yA\subseteq FA^{-1}A\subseteq F$,
which gives $y\in (F:A)$,  and so  $FA^{-1}\subseteq (F:A)$,  which leads to $
(FA^{-1})^{\ast }\subseteq (F:A)^{\ast }$.  Now we have shown that $
(F:A)^{\ast }\subseteq (F^{\ast }:A)\subseteq (FA^{-1})^{\ast }\subseteq
(F:A)^{\ast }$,  which establishes the equalities.  

(i)$\Rightarrow$(vi). 
By the proof of (i)$\Rightarrow$(v), $(F:A)^{v} \subseteq
(F^{\ast }:A)^v = (F^{\ast }:A)\subseteq (FA^{-1})^{\ast }\subseteq (F:A)^{\ast }\subseteq
(F:A)^{v}$.  This gives the required equations.

(i)$\Rightarrow$(vii).  If we insert $A^{-1}$ for $A$ in   (vi), then   we get (vii).

Next we will show that each of the conditions (iii), (iv), (v), (vi) and (vii)  implies that $F$ is $\ast$-invertible for  all   $F \in \f(D)$.

(iii)$\Rightarrow$(i). In  $ (A^{\ast }:F)=(AF^{-1})^{\ast }$,  set $A=F$ for $F \in \f(D)$.   We have $ (F^{\ast} :F)  
=(FF^{-1})^{\ast }$. Now note that $D\subseteq (F^{\ast}:F)
=(FF^{-1})^{\ast }\subseteq D$.

(iv)$\Rightarrow$(i). In  $(A^{\ast }:F^{-1}) = (AF)^{\ast }$,  set $A=F^{-1}$ for $F \in \f(D)$.

(v)$\Rightarrow$(i).  In  
$(F^{\ast }:A)=(FA^{-1})^{\ast }$,  set $A=F$ for $F \in \f(D)$.

(vi)$\Rightarrow$(i).  In  $(F^{v}:A)=(FA^{-1})^{\ast }$,  set $A=F$ for $F \in \f(D)$.     

 (vii)$\Rightarrow$(i). In  $(F^{v}:A^{-1})=(FA^{v})^{\ast }$,  set $A=F^{-1}$  for $F \in \f(D)$.

 (iii)$\Rightarrow$(viii). Applying   (iii),    we  have 
 $$
\begin{array}{rl}
((A+B):F)^{\ast
}= & ((A+B)F^{-1})^{\ast } \\
=& ((AF^{-1})^{\ast }+(BF^{-1})^{\ast })^{\ast }
=((A:F)^{\ast }+(B:F)^{\ast })^{\ast }\\
=& ((A:F)+(B:F))^{\ast }.
\end{array}
$$

 (viii)$\Rightarrow$(i).  Let  $0 \neq a, b \in D$.   Set
 $A:=(a)$,  $B:=(b)$,  and $F:=(a,b)$.  Then 
 $$
\begin{array}{rl} 
D\subseteq&  \left((a,b):(a,b)\right)^{\ast }=
\left(((a):(a,b))+((b):(a,b))\right)^{\ast } \\
=& \left(\left((a):(a,b)\right)^{\ast }+\left((b):(a,b)\right)^{\ast}\right)^{\ast}\\
=& \left(\left(a(a,b)^{-1}\right)^\ast  +\left(b(a,b)^{-1}\right)^{\ast} \right)^\ast \\
=& \left(a(a,b)^{-1} +\ b(a,b)^{-1} \right)^\ast \\
=& \left((a,b)(a,b)^{-1}\right)^{\ast} \\ \subseteq & D
\end{array}
$$
which forces $((a,b)(a,b)^{-1})^{\ast }=D$.  The conclusion follows from Theorem \ref{*prufer} ((i)$\Leftrightarrow$(ii)).

 (i)$\Rightarrow$(ix).   For all $A \in \F(D)$ and  $F, G \in \f(D)$, note that $((A:F) + (A:G))^\ast = ((A:F)^\ast + (A:G)^\ast)^\ast$; moreover $(A:F)^\ast = (AF^{-1})^\ast$ and $(A:G)^\ast = (AG^{-1})^\ast$ by (i)$\Rightarrow$(iii).  Therefore  $((A:F) + (A:G))^\ast =  ((AF^{-1})^\ast + (AG^{-1})^\ast)^\ast =   (AF^{-1} + AG^{-1})^\ast =  (A(F^{-1} + G^{-1}))^\ast  $.

Since $D$ is a $\ast$--Pr\"ufer domain, $F$ and $G$ are  $\ast$-invertible,   and thus   $F^{-1}$ and 
$G^{-1}$ are also  $\ast$-invertible. Therefore, $F\cap G$ and $F^{-1}+G^{-1}$ are $\ast$-invertible by Theorem \ref{*prufer} ((i)$\Rightarrow$(v)).  Moreover, since a $\ast$-invertible  $\ast$-ideal is a $v$-invertible $v$-ideal \cite[Proposition 3.1]{DFA},  for  $F, G \in \f^\ast(D)$,   we have in particular $F =F^\ast = F^v$, 
 $G= G^\ast = G^v$,  and  $(F^{-1}+G^{-1})^\ast = (F^{-1}+G^{-1})^v$. 
On the other hand,
$(F^{-1}+G^{-1})^{-1}  = (D: (F^{-1}+G^{-1})) = (D:F^{-1}) \cap (D: G^{-1}) = F^v \cap G^v=  F \cap G$.  
Therefore,  $(A(F^{-1} + G^{-1}))^\ast  = (A(F^{-1} + G^{-1})^\ast)^\ast = (A(F^{-1} + G^{-1})^v)^\ast  = (A((F^{-1} + G^{-1})^{-1})^{-1})^\ast  =(A(F\cap G)^{-1})^\ast $. Since $F\cap G$ is $\ast$-invertible, by Proposition \ref{*-invertible} we have $(A(F\cap G)^{-1})^\ast  = (A:(F\cap G))^\ast $. Then, putting  it  all together,  we conclude  that  (ix) holds, i.e., $((A:F) + (A:G))^\ast =(A:(F\cap G))^\ast $.

(ix)$\Rightarrow$(x). Let  $0 \neq a,b \in D$,   and set $A := (a) \cap (b), \ F := (a), \ G := (b)$. By  assumption,   we have 
$$
(D \subseteq) \ \left(((a) \cap (b)): ((a) \cap (b))\right)^\ast    =\left(\left(((a) \cap (b)):(a)\right) + \left(((a) \cap (b)):(b)\right)\right)^\ast. 
$$
  On the other hand,  
$$ \begin{array}{rl}
  \left(\left(((a) \cap (b)):(a)\right) + \left(((a) \cap (b)):(b)\right)\right)^\ast  
  =& \hskip -6pt  \left((({(a)\cap (b)}){a^{-1}})+ (({(a)\cap (b)}){b^{-1}})\right)^\ast   \\
 =& \hskip -6pt\left(((b):_D(a))+ ((a):_{D}(b))\right)^{\ast } \subseteq D \,. 
 
 \end{array}
$$
 Therefore we conclude that (x) holds.

 (x)$\Rightarrow$(i).   Note that  
 $$
 \begin{array}{rl}
 \left((a,b)(a,b)^{-1}\right)^{\ast }=& \hskip -6pt\left(a(a,b)^{-1}+b(a,b)^{-1}\right)^{\ast} \\
= & \hskip -6pt  \left(a(({(a)\cap (b)}){(ab)^{-1}})+b(({(a)\cap (b)}){(ab)^{-1}})\right)^{\ast } 
 \\
 = & \hskip -6pt  \left((({(a)\cap (b)}){b^{-1}})+ (({(a)\cap (b)}){a^{-1}})\right)^\ast  \\
 =& \hskip -6pt \left(((a):_{D}(b))+ ((b):_D(a))\right)^{\ast } 
\end{array}
$$
and apply Theorem \ref{*prufer}  ((ii)$\Rightarrow$(i)). 
\end{proof}

\begin{remark}   \bf (1) \rm  Note that from the proof of Theorem \ref{*prufer-2} it follows easily that the following conditions are equivalent to each of the conditions  (i)--(x). 
\begin{enumerate} \it

\item[(iii$^\prime$)] $(A:F)^{\ast }= (AF^{-1})^{\ast }$  for all  $A\in \F(D) $ and $F\in \f(D)$.

\item[(iii$^{\prime\prime}$)] $(A^{\ast }:F) = (AF^{-1})^{\ast }$  for all  $A\in \F(D) $ and $F\in \f(D)$.

\item[(iv$^\prime$)]  $(A:F^{-1})^{\ast }= (AF)^{\ast }$ for all  $A\in \F(D) $ and $F\in \f(D)$.

\item[(iv$^{\prime\prime}$)]  $(A^{\ast }:F^{-1}) = (AF)^{\ast }$ for all  $A\in \F(D) $ and $F\in \f(D)$.

\item[(v$^\prime$)]  $(F:A)^{\ast }=(FA^{-1})^{\ast }$ for all  $A\in \F(D) $ and $F\in \f(D)$.

\item[(v$^{\prime\prime}$)]  $(F^{\ast }:A)=(FA^{-1})^{\ast }$ for all  $A\in \F(D) $ and $F\in \f(D)$.

\item[(vi$^\prime$)]  $(F:A)^{v}= (FA^{-1})^{\ast }$ for all  $A\in \F(D) $ and $F\in \f(D)$.

\item[(vi$^{\prime\prime}$)]  $ (F^{v}:A)=(FA^{-1})^{\ast }$ for all  $A\in \F(D) $ and $F\in \f(D)$.
\end{enumerate}
\rm

As a by-product, we obtain a direct proof of \cite[Corollary 4.3]{AMZ}.

\bf (2) \rm In   analogy  with the equivalence (i)$\Leftrightarrow$(ii) of Theorem \ref{*prufer-2}, it is straightforward to prove  the following ``multiplication-type'' characterizations of the  ``Pr\"ufer-like''   classes of integral domains introduced above.

\it   
\begin{enumerate}
\item[(a)] The following properties are  equivalent.  

{\begin{enumerate}
\item[(i)]  $D$ is a $\ast$--Pr\"ufer domain. 
\item[(ii)]  If $F, G \in \f(D)$ and $F^\ast \subseteq G^\ast$, then $F^\ast = (GB)^\ast$ for some $B \in \F(D)$.
\end{enumerate}}
\item[(b)] The following properties are  equivalent.  

{\begin{enumerate}
\item[(j)]  $D$ is a $(\ast, v)$--Pr\"ufer domain. 
\item[(jj)]  If $F, G \in \f(D)$ and $F^\ast \subseteq G^v$, then $F^\ast = (GB)^\ast$ for some $B \in \F(D)$.
\end{enumerate}}
\item[(c)] The following properties are  equivalent. 
{\begin{enumerate}
\item[(i$_{\f}$)]  $D$ is a $\astf$--Pr\"ufer domain \ (= P$\ast$MD). 
\item[(ii$_{\f}$)]  If $F, G \in \f(D)$ and $F^\ast \subseteq G^\ast$, then $F^\ast = (GH)^\ast$ for some $H \in \f(D)$.
\end{enumerate}}
\item[(d)] The following properties are  equivalent. 
{\begin{enumerate}
\item[(j$_{\f}$)]  $D$ is a $(\astf, v)$--Pr\"ufer domain. 
\item[(jj$_{\f}$)]  If $F, G \in \f(D)$ and $F^\ast \subseteq G^v$, then $F^\ast = (GH)^\ast$ for some $H \in \f(D)$.
\end{enumerate}}
\end{enumerate}

\end{remark}

\begin{remark} \rm Referring  to Theorem 6.6 in \cite{LM}, which provides several characterizations of Pr\"ufer domains, we can summarize that conditions (2), (5), (6), (7),  (8),   and (9)  of that theorem have been modified in a canonical way (see, respectively,  conditions (ii), \ (iv$_{\f}$)\&(iv$_{\!\f\!\F}$), \ (iii$_{\f}$)\&(iii$_{\F}$) \ of Theorem  \ref{*prufer} and conditions \ (ii),  \ 
 \ (viii),   \  (ix)  \   of Theorem \ref{*prufer-2}) in order to obtain characterizations of   $\ast$--Pr\"ufer domains.      

We have also observed  that  condition (3)  (of  \cite[Theorem 6.6]{LM})   extends  in a natural way   to   ``$\ast$ is an a.b. star operation'',   and  we have just seen that  $D$ being $\ast$--Pr\"ufer   implies  that $\ast$ is an a.b. star operation on $D$,  but not conversely (Proposition \ref{*prufer-necessary} (3) and Remark \ref{rk:2.5} (c)). 

 Moreover,   there is no natural  modification of condition (4) (of  \cite[Theorem 6.6]{LM}) which can  provide a characterization of    $\ast$--Pr\"ufer domains:    take, for instance, a  one-dimensional  quasi-local CICD  (hence,  $v$--domain)   which is not a valuation domain  \cite{N-1, N-2, R}.  We were unable to find an appropriate modification of  condition
(10) (of  Theorem 6.6 in \cite{LM}) leading to a characterization of $\ast$--Pr\"ufer domains.   
\end{remark}

\begin{theorem} \label{*-inv} Let $D$ be an integral domain with quotient field $K$, let  $\ast$ be a star operation on  $D$,   and let $\Inv^{\ast}(D)$ be the group of $\ast$-invertible $\ast$-ideals of $D$  under  $\ast$-multiplication (defined by $A\ast\! B:= (AB)^\ast$  for  all $A, B \in \Inv^{\ast}(D)$). 
\begin{enumerate}
\item[(a)] $D$ is a
  $\ast$--Pr\"ufer domain   if and only if $\Inv^{\ast}(D)$  is a  lattice-ordered   abelian   group  under the relation $A\leq B$ 
defined by $A\supseteq B$   for  $A, B \in \Inv^{\ast}(D)$   with     \ $\sup (A,B)=A\cap B $ \  and \ $\inf (A,B)=(A+B)^{\ast}=(A+B)^{v}$ \ for   all $A,B\in \Inv^{\ast}(D)$.
\item[(b)] If $\Inv^{\ast}(D)$  is a  lattice-ordered   abelian   group (under the relation $\leq $ 
defined  above),   then   $\sup (A,B)=A\cap B $ and $\inf (A,B)=(A+B)^{v}$  for  all $A,B\in \Inv^{\ast}(D)$. In this situation, $D$ is a   $(\ast, v)$--Pr\"ufer domain.  

 \item[(c)] $D$ is a
  $\ast$--Pr\"ufer domain  if and only if $D$ is a  $(\ast, v)$--Pr\"ufer domain   and $(A+B)^{\ast}=(A+B)^{v}$  for all $A,B\in \Inv^{\ast}(D)$.  

\end{enumerate}

\end{theorem} 

\begin{proof} That $\Inv^{\ast}(D)$ is   an abelian  group was observed in \cite[page 812]{DFA}. That $\Inv^{\ast}(D)$
is a partially ordered group (under the partial ordered defined above) is easy to see because for $A,B\in \Inv^{\ast}(D)$ and
for any nonzero fractional ideal $X$, $A\supseteq B$ implies $XA\supseteq XB$
and hence $(XA)^{\ast}\supseteq (XB)^{\ast}.$ Thus, in particular,  for all $X,A,B\in \Inv^{\ast}(D)$, $A\leq B$ implies $X\!\ast\! A = (XA)^{\ast}\leq (XB)^{\ast} = X \!\ast\! B$.  So  the relation $\leq $ is compatible with
group  multiplication,   and hence  $\Inv^{\ast}(D)$ is a partially ordered group \cite[pages 61 and 107]{FS}.

 (a)  Now suppose that $D$ is a    $\ast$--Pr\"ufer domain.   By  Theorem \ref{*prufer}  ((i)$\Rightarrow$(iii$_{\F}$)),   
$A\cap B \ (= (A\cap B)^{\ast})$ and $(A+B)^{\ast}$ both belong to $\Inv^{\ast}(D),$ whenever $A,B\in
\Inv^{\ast}(D)$.
Therefore, it is straightforward to verify that   $A\cap B =\sup (A,B)$ and $(A+B)^{\ast}=\inf (A,B)$  for $A,B\in \Inv^{\ast}(D)$.
 Thus $\Inv^{\ast}(D)$ is
 a  lattice-ordered  group \cite[page 107]{FS}.  Note also that, since a $\ast$-invertible $\ast$-ideal is a $v$-ideal \cite[Corollaire 1, page 24]{J}, $(A+B)^{\ast} = ((A+B)^{\ast})^v = (A+B)^v$. Therefore
 $\inf (A,B)= (A+B)^{\ast}=(A+B)^v$   for  $A,B\in \Inv^{\ast}(D)$.

 Conversely,   suppose that $\Inv^{\ast}(D)$ is a  lattice-ordered  group (under $\leq $
defined above) and that $\inf (A,B) =(A+B)^\ast$   for   $A,B\in \Inv^{\ast}(D)$.
In particular,  $\inf (aD, bD)$  $ =(aD+bD)^{\ast} \in \Inv^{\ast}(D)$   for all  $0\neq a,b\in
D$;  hence
every two generated nonzero ideal is $\ast$-invertible,   and thus  $D$ is    a $\ast$--Pr\"ufer domain  by 
Theorem \ref{*prufer} ((ii)$\Rightarrow$(i)).

 (b)   We start  by showing that under the present assumption  $\inf (A,B)= (A+B)^v$ 
for $A,B\in \Inv^{\ast}(D)$. Since  $\inf (A,B) \in \Inv^{\ast}(D)$  and, as we observed above,  every   $\ast$-invertible $\ast$-ideal is a $v$-ideal,   clearly   $\inf (A,B) \supseteq (A+B)^v$. For the reverse containment, for  all   $H \in \Inv^{\ast}(D)$ such that $H \supseteq A$ and $H\supseteq B$, we have that $H \supseteq \inf(A,B)$.  Since $\Inv^{\ast}(D)$ contains all principal
fractional ideals,  in particular   we have $\bigcap \{zD \mid0 \neq z \in K, \ zD \supseteq A \mbox{ and } zD\supseteq B \} \supseteq \bigcap \{H \in \Inv^{\ast}(D)\mid H \supseteq A \mbox{ and } H\supseteq B \}  \supseteq \inf (A,B)$.  Since  $\bigcap \{zD \mid0 \neq z \in K, \ zD \supseteq A \mbox{ and } zD\supseteq B \} = (A+B)^v$, we conclude that $\inf (A,B) = (A+B)^v$.

Next we show that  $\sup (A,B)= A\cap B$   for  all $A, B\in \Inv^{\ast}(D)$. It is easy to verify that 
$(\sup(A,B)A^{-1}B^{-1})^\ast  =\sup(A^{-1},B^{-1})$ and that $\sup(A^{-1},B^{-1})=(\inf(A, B))^{-1}$ since,   for all $A, B\in \Inv^{\ast}(D)$,   $A^{-1}, B^{-1} \in\Inv^{\ast}(D)$ and  $A \leq B$ if and only if $A^{-1} \geq B^{-1}$. Therefore we have 
$(\sup(A,B)\inf(A,B))^\ast = (AB)^\ast$ or, equivalently,  $(\sup(A,B)\inf(A,B)A^{-1}B^{-1})^\ast = D$  for  all $A, B\in \Inv^{\ast}(D)$. 
Replacing $\inf(A,B)$ by $(A+B)^v$ and applying the $v$-operation on both sides, we have ~$ (\sup(A,B)(A+B)^vA^{-1}B^{-1})^v = ((\sup(A,B)(A+B)^vA^{-1}B^{-1})^\ast)^v =D$.  So
$D = (\sup(A,B)(A+B)^vA^{-1}B^{-1})^v  = (\sup(A,B)(A+B)A^{-1}B^{-1})^v = (\sup(A,B)$ $(B^{-1}+A^{-1}))^v $, forcing $\sup(A,B) = (B^{-1}+A^{-1})^{-1} = A^v \cap B^v = A \cap B$ (since $A$ and $B$ are $v$-ideals, as observed above).

In order to show that $D$ is a  $(\ast, v)$--Pr\"ufer domain,   we start by showing that if $F$ is a nonzero two generated ideal of  $D$,   then $F^v$ is $\ast$-invertible. Let $F := aD +bD$, with $a, b \in D$ and  $ab \neq 0$.  We know that  
$abD = (\inf(aD, bD)\sup(aD, bD))^\ast = ((aD \cap bD)(aD+bD)^v)^\ast$,   and thus  $D = ((b^{-1}D \cap a^{-1}D)(aD+bD)^v)^\ast$, i.e., $(aD+bD)^v$ is $\ast$-invertible.
The general case can be obtained by induction. Let $F$ be a nonzero ideal of $D$ generated by $n\geq 2$ elements and let $c \in D \setminus F$. By the previous  arguments,   we have that the ideal $F +cD$, generated by $(n+1)$ elements,  is such that $ (F +cD)^v  = \inf(F, cD)  \in \Inv^{\ast}(D)$.

 (c)  The ``only if part'' follows immediately from (a) and (b). For the ``if part'', let $D$ be a $(\ast, v)$--Pr\"ufer domain and let $F:= aD + bD$, where  $0\neq a, b \in D$.   Since $aD, bD \in \Inv^\ast(D)$, by assumption $F^v = ( aD + bD)^v = ( aD + bD)^\ast = F^\ast$. Therefore $D = (F^vF^{-1})^\ast = (F^\ast F^{-1})^\ast  = (FF^{-1})^\ast $. The conclusion is an immediate consequence of Theorem \ref{*prufer} ((ii)$\Rightarrow$(i)). 
\end{proof}

With the proof of   Theorem \ref{*-inv},    we have amply established
the existence of a sort of GCD   in $\Inv^{\ast }(D)$  for   each pair of elements of $\Inv^{\ast }(D)$  
when  $D$ is a $\ast $--Pr\"ufer domain.  However, the results
are in terms of \ inf \ and \ sup \ of elements of the  lattice-ordered  group $ 
\Inv^{\ast }(D).$ We now establish the existence of  the  ($\ast $-invertible $\ast 
$-ideal) GCD of $\ast $-invertible integral $\ast $-ideals   using  purely ring-theoretic   means. 

Before we do that,  let us note that old masters such as van
der Waerden regarded an integral ideal $A$ of an integral domain $D$  as   a
\it divisor \rm of another ideal $B$ of $D$ if $A\supseteq B$, extending the well known  property that, for  $0 \neq a, b \in D$,   $aD \supseteq bD$ if and only if $a | b$.   In turn, the ideal  $B$ could be
termed as a \it multiple \rm of the ideal  $A$.   Now,  given two integral ideals $A,B$ of  $D$,   the
ideal $A+B$ has the property that $A+B$ is a divisor of $A$ and $B$ and any
common divisor $C$ of $A$ and  $B$   contains $A$ and  $B$,  and hence $A+B.$ In
other  words,   any common divisor  of   $A$ and $B$ is a divisor of $A+B$. Thus 
$A+B$ fitted the bill as \it the greatest common divisor \rm  of $A$ and $B$. 
In a
similar fashion $A\cap B$ was regarded as \it the least common multiple \rm  of $A$ and $B$ \cite[
Vol. 2, page 119]{Waerden}. 

Now the trouble with this approach is that it is too
general and so can only work in a very strict  environment   such as a PID or a
Dedekind domain, the kind of rings the \textquotedblleft
ancients\textquotedblright\ worked with.  Besides,   there were other ways of
looking at GCD's, such as generalizations of the GCD of two integers, which
is an integer.   Also,   if we are dealing with $\boldsymbol{\mathcal{I}}^\ast \!(D)$, the set of
integral $\ast $-ideals of $D$,  and we want the GCD of two $\ast $-ideals $A, B\in \boldsymbol{\mathcal{I}}^\ast \!(D)$  to belong  to $\boldsymbol{\mathcal{I}}^\ast \!(D)$, then in general  $A+B$ would not deliver the \textquotedblleft greatest common
divisor\textquotedblright  \   in   $\boldsymbol{\mathcal{I}}^\ast (D)$,  in the language of van der
Waerden.  So to find the GCD of $A,B\in
\boldsymbol{\mathcal{I}}^\ast \!(D)$ inside  $\boldsymbol{\mathcal{I}}^\ast \!(D)$,   we need to consider $(A+B)^{\ast }$,  which may be a
proper divisor of $A+B$.  In other words,  we need GCD's from a pre-assigned
set.  Of  course,   we also need our GCD to  be   something like the GCD in  Pr\"{u}fer domains  that we defined in the introduction. Having established
what we want,  we state  a GCD-type characterization of $\ast$--Pr\"ufer domains.

\begin{proposition} Let $D$ be an integral domain,  $\ast$  a star operation on  $D$,   and let $\Inv_{\boldsymbol{\mathcal{I}}}^{\ast }(D)$ be the set of integral $\ast $-invertible $\ast $-ideals of $D$. 

Assume that  $D$ is a $\ast $--Pr\"{u}fer domain.  If   $A, B\in \Inv_{\boldsymbol{\mathcal{I}}}^{\ast }(D)$,   then there is a unique ideal $
C:=(A+B)^{\ast }\in \Inv_{\boldsymbol{\mathcal{I}}}^{\ast }(D)$ such that $A=(A_{1}C)^{\ast }$, $
B=(B_{1}C)^{\ast }$, where $(A_{1}+B_{1})^{\ast }=D$.
\ Conversely, if $D$ is an
integral domain such that for  all  $A,B\in \Inv_{\boldsymbol{\mathcal{I}}}^{\ast }(D)$,   there is
a unique ideal $C\in \Inv_{\boldsymbol{\mathcal{I}}}^{\ast }(D)$ such that $A=(A_{1}C)^{\ast }$, $
B=(B_{1}C)^{\ast }$, where $(A_{1}+B_{1})^{\ast }=D$, then $D$ is a $
\ast $--Prufer domain.
\end{proposition}

\begin{proof} Let us first note that if $I$ is an integral $\ast $-invertible $\ast 
$-ideal of a $\ast $--Pr\"{u}fer domain $D$ and $J$ is an ideal contained in  $
I$,   then $J^{\ast }=(IH)^{\ast }$ for some integral ideal $H$ of $D$.  This follows
since $J\subseteq I$ implies $JI^{-1}=:H \subseteq D$.   Now, multiplying both
sides of the equality  $JI^{-1}=H$ by $I$ and  applying $\ast$,  we get $J^{\ast
}=(IH)^{\ast }$.

 Next, let $C:=(A+B)^{\ast }$.  Since $D$ is a $\ast 
$--Pr\"{u}fer domain,  $C$ is $\ast $-invertible (Theorem \ref{*prufer}  ((i)$\Rightarrow$(vi)). Now as  $A,B \subseteq C$,  we have $
(C^{-1}A)^{\ast }=: A_{1}\subseteq D$ and $(C^{-1}B)^{\ast }=:B_{1}\subseteq D$ so that $
A= A^\ast= (CA_{1})^{\ast }$ and $B=B^\ast=(CB_{1})^{\ast }$.   Now $C=((CA_{1})^{\ast
}+(CB_{1})^{\ast })^{\ast }=(C(A_{1}+B_{1}))^{\ast }$. Multiplying both
sides of the equality $C=(C(A_{1}+B_{1}))^{\ast }$ by $C^{-1}$  and   applying $\ast $,  we
get $(A_{1}+B_{1})^{\ast }=D$. 

The proof of the converse entails showing
that  for all  $A,B\in  \Inv_{\boldsymbol{\mathcal{I}}}^{\ast }(D)$, $ (A+B)^{\ast }\in
 \Inv_{\boldsymbol{\mathcal{I}}}^{\ast }(D)$ (Theorem \ref{*prufer} (ii)$\Rightarrow$(i))).  By  assumption,   we have   $ (A+B)^{\ast }= ((CA_{1})^{\ast
}+(CB_{1})^{\ast })^{\ast }=(C(A_{1}+B_{1})^{\ast })^{\ast } =C^\ast = C \in \Inv_{\boldsymbol{\mathcal{I}}}^{\ast }(D)$.  \end{proof}

 As a  consequence of   Theorem \ref{*-inv},  we  have  

\begin{corollary} \label{cor:*-inv} Let $D$ be an integral domain.  
\begin{enumerate}
\item [(a)] 
$D$ is a $v$--domain (respectively, a P$v$MD, a  generalized  GCD domain) if and only if \  $\Inv^v(D)$ (respectively, $\Inv^t(D)$, $\Inv(D)$)  is a  lattice-ordered   abelian   group (under $\leq $
defined in  Theorem \ref{*-inv}).
\item [(b)]  Assume  that $\Inv(D)$    is  lattice-ordered  and that $A+B = (A+B)^v$ for all $A, B \in \Inv(D)$.  Then $D$ is a Pr\"ufer domain and, clearly,  $\Inv(D) = \f(D)$.  
\end{enumerate}
\end{corollary}
 \begin{proof} (a) The ``if part'' is  a    consequence of  Theorem \ref{*-inv} (b) and Example \ref{ex:2.1}. The ``only if part'' for $v$--domains  and  P$v$MD's follows from Theorem \ref{*-inv} (a) (and from Example \ref{ex:2.1}). If $D$ is a  generalized  GCD domain ( = $(d, v)$--Pr\"ufer domain), then   it is well known that $\Inv(D)$ is  lattice-ordered  and $\Inv(D) = \{F^v \mid F \in \f(D) \}$   \cite[Theorem 1 ((1)$\Rightarrow$(4))]{AA}.

(b) is an easy consequence of Theorem \ref{*-inv} (c). 
\end{proof}

\begin{remark} 
\rm The ``P$v$MD part'' of Corollary \ref{cor:*-inv}   gives back a classical
characte\-rization of these domains (see e.g. \cite[page 55]{J}, \cite[page 717]{Gr} and 
\cite[Proposition 2.4]{Z-inv}). 
As we mentioned above, the ``GGCD part''  is well known \cite{AA}.  On the other hand,  the ``$v$--domain part'' is new. 
\end{remark}


\begin{corollary} \label{f^v-group}
 Let $D$ be an   integral domain. 
 \begin{enumerate}
 
 \item[(a)] Assume that $D$ is a CICD. Then   $D$ is a $v$--domain  with   $\F^{v}(D)=\Inv^{v}(D)$, and  moreover,  $\F^{v}(D)$  is  a complete  lattice-ordered   abelian  
group (under the order $\leq $
defined by \ $I\leq J$  $:\Leftrightarrow I\supseteq J$  \    for all $I,J \in \F^v(D)$).
 \item [(b)]  Assume that $D$ is  a  $(t,v)$--Dedekind domain   
(Example \ref{ex:1.7}). Then $D$ is  a completely integrally closed P$v$MD  with  $\F^{v}(D)=\Inv^{t}(D)$,   and  moreover, $
\F^{v}(D)$ is a complete  lattice-ordered   abelian  group (under the order $\leq $
defined above).

\item[(c)]  Assume that   $D$ is a pseudo-Dedekind domain (i.e.,  a  $(d,v)$--Dedekind domain,  
Example \ref{*=d}). Then $D$  is a   generalized    GCD domain (i.e., a   $(d,v)$--Pr\"ufer  domain,    Example \ref{ex:2.1}) with  $\F^{v}(D)=\Inv(D)$,  and  moreover,  $
\F^{v}(D)$  is a complete  lattice-ordered   abelian  group (under the order $\leq $
defined above).

\end{enumerate}
\end{corollary}

\begin{proof} 
(a) Since a CICD (= $v$--CICD, by Example \ref{ex:1.6})  is a $v$--domain, $\Inv^{v}(D)$ is a lattice-ordered group from
Corollary \ref{cor:*-inv}. For the  completeness,  recall that a
lattice-ordered group $G$ is said to be complete if every nonempty subset of 
$G$ that is bounded from    below has a  greatest  lower  bound (or, equivalently,
 if
every nonempty subset of $G$ that is bounded from   above has a least
upper  bound). Note   that   when $D$ is completely integrally closed, the
 lattice-ordered  group $\Inv^{v}(D)$ coincides with $\boldsymbol{F}
^{v}(D)$.   Let $\{A_{\lambda } \mid \lambda \in \Lambda\}$  be a nonempty
collection of ideals   in  $\Inv^{v}(D)$  bounded   below  in $\Inv^{v}(D)$,    that is,   there is  
$J\in\Inv^{v}(D)$ such that $A_{\lambda }\geq J$ for all $\lambda \in \Lambda$.  In other words, 
 $A_{\lambda }\subseteq J$  for all $\lambda \in \Lambda$. Then $\sum_\lambda A_{\lambda}\subseteq J$, and hence $\left(\sum_\lambda A_{\lambda}\right)^{v} \subseteq J^v = J$.  This gives 
 $A_{\lambda }\subseteq \left(\sum_\lambda A_{\lambda}\right)^{v} \subseteq J$, which
translates \ (in $(\Inv^{v}(D), \leq)$) \  to $A_{\lambda}\geq \left(\sum_\lambda A_{\lambda}\right)^{v} \geq J$. 
  Since $\left(\sum_\lambda A_{\lambda}\right)^{v} \in \F^{v}(D)= \Inv^{v}(D)$, we ~conclude
that $\left (\sum A_{\lambda}\right)^{v}$ is a lower bound and, more precisely, it is easy to verify 
that $\left(\sum A_{\lambda}\right)^{v}$ is in fact the greatest lower bound of $\{A_{\lambda} \mid \lambda \in \Lambda \}$.

(b) In this   case,   it is clear that $\F^v(D) = \Inv^t(D)$.  Since ``$t$-invertible'' implies  ``$v$-invertible",  we have that $\F^v(D) = \Inv^v(D) \  ( = \Inv^t(D))$,   and thus $D$ is completely integrally closed.   Moreover a $(t,v)$--CICD is a particular  $(t,v)$--Pr\"ufer domain,   which is a P$v$MD (Example \ref{ex:2.1}). We conclude by (a).

(c) The proof is similar: in this  case,  $\F^v(D) = \Inv(D)$ and hence,  in particular, $\F^v(D) = \Inv^v(D) =\Inv(D)$.
\end{proof}

 Note that  the converse of  part (c) of the previous corollary  is also true.  In   fact,   it is known that an integral domain is pseudo-Dedekind if and only if $\Inv(D)$ is a complete  lattice-ordered abelian  group \cite[Theorem 2.8]{AK}.

\begin{corollary} \label{v-inv-prin}
 Let $D$ be an integral domain with quotient field $K$. Assume that  $D$ is a pseudo-principal domain (i.e, $A^{v}$ is principal for  all   $
A\in \F(D)$). Then $D$ is a GCD domain such that  $\F^{v}(D)$ is isomorphic to    the group of divisibility  of $D$.  Moreover,  $\F^{v}(D)$ is a complete  lattice-ordered   abelian   group (under the order $\leq $ defined by $I\leq J$ if $I\supseteq J$).
\end{corollary}
\begin{proof} Recall that the group of divisibility of $D$ is the multiplicative  abelian  group $G(D):=K^{\times}/{{\mathcal U}}(D)$, where $K^{\times} := K \setminus \{0\}$ and ${{\mathcal U}}(D)$ is the group of units of $D$, endowed with   a partial order   defined by $x{{\mathcal U}}(D) \leq y{{\mathcal U}}(D)$ if $yx^{-1} \in D$.
It is easy to see that the group of divisibility of $D$ is canonically isomorphic to 
$ \Prin(D):=\{zD\mid 0\neq z \in
K\}$  with a partial order $\leq$ defined by  $xD\leq yD$ if $xD\supseteq yD$ \cite[page 172]{G}.

Since a GCD domain is characterized by the fact that   $F^v$ is principal for all $F \in \f(D)$  (cf. for instance \cite[Proposition 1.19]{Pi} or   \cite[Remark 2.2]{DDA}),   it is straightforward ~from the assumption that $D$ is a GCD domain and $\F^{v}(D)=\Prin(D)$,  with identical definitions of partial order.

Next, 
for every nonempty subset ${\boldsymbol{\mathcal S}}$ of principal
fractional ideals of $D$  bounded below under $\leq$\,,\  let $A$ be the fractional ideal of $D$   generated by  the ideals in  ${\boldsymbol{\mathcal S}}$. Then    $A^{v}$ is
principal by assumption,   and   $  A^v  \supseteq sD$ for all $sD \in {\boldsymbol{\mathcal S}}$.  Thus $A^v \leq sD$ (in $(\F^{v}(D), \leq))$  for all $sD \in {\boldsymbol{\mathcal S}}$. It is routine to show  that  the   principal fractional ideal $A^v$ is   in fact the  greatest lower  bound of the family ${\boldsymbol{\mathcal S}}$.  
\end{proof}

\end{document}